\documentclass[a4paper,10pt]{article}

\usepackage{t1enc}
\usepackage[latin1]{inputenc}
\usepackage{graphics}
\usepackage{graphicx}
\usepackage{amsfonts}
\usepackage{eufrak}
\usepackage{xypic}
\usepackage{amssymb}
\usepackage{amsmath}
\usepackage{amsthm}
\usepackage{mathrsfs}

\theoremstyle{plain}
\newtheorem{thm}{Theorem}[section]
\newtheorem{lemma}[thm]{Lemma}
\newtheorem{prop}[thm]{Proposition}
\newtheorem{cor}[thm]{Corollary}

\theoremstyle{definition}
\newtheorem{df}[thm]{Definition}
\newtheorem{rem}[thm]{Remark}
\newtheorem{ex}[thm]{Example}

\renewcommand{\a}{\mathcal{A}}
\renewcommand{\b}{\mathcal{B}}

\providecommand{\gm}{\Gamma}
\providecommand{\h}{\mathcal{H}}

\title{\bf On Enveloping $C^*$-Algebras of Hecke Algebras}
\author{Rui Palma\\ {\footnotesize \emph{Department of Mathematics, University of Oslo,}}\\  {\footnotesize \emph{P.O. Box 1053 Blindern, NO-0316 Oslo, Norway}}\\ {\footnotesize\emph{E-mail: ruip@math.uio.no}}}
\date{}

\begin{document}

\maketitle

\begin{abstract}
We give a sufficient condition for a $^*$-algebra with a specified basis to have an enveloping $C^*$-algebra. Particularizing to the setting of a Hecke algebra $\h(G, \gm)$, we show that under a suitable assumption not only we can assure that an enveloping $C^*$-algebra $C^*(G, \gm)$ exists, but also that it coincides with $C^*(L^1(G,\gm))$, the enveloping $C^*$-algebra of the $L^1$-Hecke algebra. Our methods are used to show the existence of $C^*(G, \gm)$ and isomorphism with $C^*(L^1(G,\gm))$ for several classes of Hecke algebras. Most of the classes which are known to satisfy these  properties are covered by this approach, and we also describe some new ones.
\end{abstract}

{\renewcommand{\thefootnote}{}
\footnotetext{\emph{Date:} \today}}

{\renewcommand{\thefootnote}{}
\footnotetext{Research supported by the Research Council of Norway, the Nordforsk research network ``Operator Algebra and Dynamics'' and Funda\c c\~ ao para a Ci\^ encia e Tecnologia grant SFRH/BD/43276/2008.}}

\section*{Introduction}

A \emph{Hecke pair} $(G, \gm)$ consists of a group $G$ and a subgroup $\gm \subseteq G$ for which every double coset $\gm g \gm$ is the union of finitely many left cosets. In this case $\gm$ is also said to be a \emph{Hecke subgroup} of $G$. Examples of Hecke subgroups include finite subgroups, finite-index subgroups and normal subgroups. Hecke subgroups are also sometimes called \emph{almost normal subgroups} (although we will not use this terminology here) and it is in fact many times insightful to think of this definition as a generalization of the notion of normality of a subgroup.

Given a Hecke pair $(G,\gm)$ the \emph{Hecke algebra} $\h(G, \gm)$ is a $^*$-algebra of functions over the set of double cosets $\gm \backslash G / \gm$, with a suitable convolution product and involution. It generalizes the definition of the group algebra $\mathbb{C}(G / \gm)$ of the quotient group when $\gm$ is a normal subgroup.

 The interest in Hecke algebras in the realm of operator algebras was to a large extent raised through the work of Bost and Connes \cite{bost connes} on phase transitions in number theory, and since then several authors have studied $C^*$-algebras which arise as completions of Hecke algebras. There are several canonical $C^*$-completions of a Hecke algebra (\cite{tzanev}, \cite{schl}) and the question of existence of a maximal one (i.e. an enveloping $C^*$-algebra) has been of particular interest (\cite{bost connes}, \cite{arledge}, \cite{brenken}, \cite{hall}, \cite{tzanev}, \cite{glockner}, \cite{laca larsen}, \cite{compac}, \cite{schl}, \cite{exel}). One of the reasons for that, firstly explored by Hall \cite{hall}, has to do with how well $^*$-representations of a Hecke algebra $\h(G, \gm)$ correspond to unitary representations of the group $G$ that are generated by their $\gm$-fixed vectors. It was shown by Hall \cite{hall} that for such a correspondence to hold it is necessary that the Hecke algebra has an enveloping $C^*$-algebra, which does not always happen. It was later clarified by Kaliszewski, Landstad and Quigg \cite{schl} that such a correspondence holds precisely when an enveloping $C^*$-algebra exists and coincides with other canonical $C^*$-completions.

The problem of deciding if a Hecke algebra has an enveloping $C^*$-algebra seems to be of a non-trivial nature, with satisfactory answers, arising from various distinct methods, known only for certain classes of Hecke pairs.

The main motivation for the present article is to give a unified approach to this problem for a large class of Hecke pairs. We recover most of the known cases in the literature but also several new ones. We achieve this by associating a directed graph to a Hecke algebra $\h(G, \gm)$, whose vertices are the double cosets and whose directed edges are determined by how products of the form $(\gm g \gm)^* *\gm g \gm$ decompose as sums of double cosets. We prove that finiteness of the \emph{co-hereditary set} generated by a vertex $\gm g \gm$, i.e. the set of vertices one encounters by moving forward in the graph starting from $\gm g \gm$, implies that
\begin{align*}
 \sup_{\pi} \| \pi(\gm g \gm)\| < \infty\,,
\end{align*}
where the supremum runs over the $^*$-representations of the Hecke algebra. Thus, analysing these co-hereditary sets gives valuable information regarding the existence of enveloping $C^*$-algebras. Moreover, we prove that if \emph{all} double cosets generate finite co-hereditary sets, then the enveloping $C^*$-algebra of $\h(G, \gm)$ exists and coincides with $C^*(L^1(G, \gm))$, the enveloping $C^*$-algebra of the $L^1$-Hecke algebra (one of the canonical $C^*$-completions).

We develop certain tools, based on iterated commutators on the group $G$, that allow us to show that our assumptions hold in a variety of classes of Hecke pairs, and thus enable us to answer affirmatively the question of existence of enveloping $C^*$-algebras of the corresponding Hecke algebras. Some of the new results we prove state that if a group $G$ satisfies some generalized nilpotency property, then for $\emph{any}$ Hecke subgroup $\gm$ the Hecke algebra $\h(G, \gm)$ has an enveloping $C^*$-algebra which coincides with $C^*(L^1(G, \gm))$. These results will enable us to show, in a following article (see \cite{palma2}), that for any $G$ satisfying such properties, Hall's correspondence holds for any Hecke subgroup.

We also notice that the classes of Hecke algebras studied in the present work, and therefore most of those studied in the literature, satisfy a stronger property than just having an enveloping $C^*$-algebra: they are in fact $BG^*$-algebras. The standard reference for this class of $^*$-algebras is Palmer \cite{palmer}, but we also give a short description in Section 1. The reason for considering this stronger property is not only because of how well-behaved these $^*$-algebras are, but also because it is natural to consider $BG^*$-Hecke algebras in the context of crossed products by Hecke pairs (see \cite{palmathesis}).

The paper is organized as follows: in Section \ref{section preliminaries} we set up the conventions, notation, and background results regarding $^*$-algebras, Hecke algebras and also directed graphs, that will be used throughout the article. In the setting of directed graphs the most important notion will be that of a \emph{co-hereditary set}.

In Section \ref{graph assoc with star alg} we associate a directed graph to any $^*$-algebra with a given basis and prove our first main result, which states that we can put a bound on the norm of all representations of the elements of a finite co-hereditary set. As a consequence, if the $^*$-algebra with the given basis is generated by its finite co-hereditary sets, it must have an enveloping $C^*$-algebra.

In Section \ref{hecke star alg section} we prove the second main result of this article: that in the case of a Hecke algebra $\h(G, \gm)$ if all double cosets generate finite co-hereditary sets, then an enveloping $C^*$-algebra exists and it coincides with $C^*(L^1(G, \gm))$.

In Section \ref{methods section} we present some tools for determining if a co-hereditary set generated by a double coset is finite. These methods will be then used in Section \ref{classes of Hecke pairs} to study the existence of an enveloping $C^*$-algebra, and the isomorphism with $C^*(L^1(G, \gm))$, for several classes of Hecke pairs $(G, \gm)$.

In Section \ref{reduction section} we show that the problem of existence of an enveloping $C^*$-algebra for $\h(G, \gm)$ can be reduced to the same problem but for a smaller Hecke subalgebra $\h(H,\gm)$, where $\gm \subseteq H \subseteq G$ is an ascendant subgroup.

Finally in Section \ref{questions section} we give some concluding remarks and state some open questions.

The present work is part of the author's Ph.D. thesis \cite{palmathesis} written at the University of Oslo. The author would like to thank his advisor Nadia Larsen for the very helpful discussions, suggestions and comments during the elaboration of this work.\\

\section{Preliminaries}
\label{section preliminaries}

\subsection{Preliminaries on $^*$-Algebras}

 Let $\mathscr{V}$ be an inner product space over $\mathbb{C}$. Recall that a function $T: \mathscr{V} \to \mathscr{V}$ is said to be \emph{adjointable} if there exists a function $T^*: \mathscr{V} \to \mathscr{V}$ such that
\begin{align*}
 \langle T\xi \,,\, \eta \rangle = \langle \xi \,,\, T^* \eta \rangle\,,
\end{align*}
for all $\xi, \eta \in \mathscr{V}$. Recall also that every adjointable operator $T$ is necessarily linear and that $T^*$ is unique and adjointable with $T^{**} = T$. We will use the following notation:
\begin{itemize}
 \item $L(\mathscr{V})$ denotes the $^*$-algebra of all adjointable operators in $\mathscr{V}$
 \item $B(\mathscr{V})$ denotes the $^*$-algebra of all bounded adjointable operators in $\mathscr{V}$.
\end{itemize}
Of course, we always have $B(\mathscr{V}) \subseteq L(\mathscr{V})$, with both $^*$-algebras coinciding when $\mathscr{V}$  is a Hilbert space (see, for example, \cite[Proposition 9.1.11]{palmer}).

Following \cite[Def. 9.2.1]{palmer}, we define a \emph{pre-$^*$-representation} of a $^*$-algebra $\a$ on an inner product space $\mathscr{V}$ to be a $^*$-homomorphism $\pi: \a \to L(\mathscr{V})$ and a \emph{$^*$-representation} of $\a$ on a Hilbert space $\mathscr{H}$ to be a $^*$-homomorphism $\pi: \a \to B(\mathscr{H})$. As in \cite[Def. 4.2.1]{palmer11}, a pre-$^*$-representation $\pi:\a \to L(\mathscr{V})$ is said to be \emph{normed} if $\pi(\a) \subseteq B(\mathscr{V})$, i.e. if $\pi(a)$ is a bounded operator for all $a \in \a$. We now make a seemingly similar definition, but where the focus is on the elements of the $^*$-algebra, instead of its pre-$^*$-representations:\\

\begin{df}
 Let $\a$ be a $^*$-algebra. We will say that an element $a \in \a$ is \emph{automatically bounded} if $\pi(a) \in B(\mathscr{V})$ for any pre-$^*$-representation $\pi:\a \to L(\mathscr{V})$.\\
\end{df}

Easy examples of automatically bounded elements in a $^*$-algebra are unitaries, projections, or more generally, partial isometries.

Given a $^*$-algebra $\a$ let
\begin{align*}
 \a_b & :=  \{a \in \a : a \;\;\text{is automatically bounded} \}\,.\\
\end{align*}

\begin{df}[\cite{palmer}, Def. 10.1.17]
 A $^*$-algebra $\a$ is called a \emph{$BG^*$-algebra} if every element $a \in \a$ is automatically bounded, i.e. if $\a_b = \a$. Equivalently, $\a$ is a \emph{$BG^*$-algebra} if all pre-$^*$-representations of $\a$ are normed.\\
\end{df}

The function $\| \cdot \|_u: \a \to \mathbb{R}^+_0 \cup \{\infty\}$ defined by
\begin{align*}
\|a\|_u := \sup_{\pi} \| \pi(a)\|\,,
\end{align*}
where the supremum is taken over all $^*$-representations of $\a$, will be called the \emph{universal norm} of $\a$. An element $a \in \a$ will be said to have a \emph{bounded universal norm} if $\|a\|_u < \infty$, and the set of all elements $a \in \a$ which have a bounded universal norm will be denoted by $\a_u$, i.e.
\begin{align*}
 \a_u & := \{a \in \a : \|a\|_u < \infty \}\,.
\end{align*}

When $\a_u = \a$ the universal norm becomes a true $C^*$-seminorm, being actually the largest possible $C^*$-seminorm in $\a$. The Hausdorff completion of $\a$ in the universal norm is then a $C^*$-algebra called the \emph{enveloping $C^*$-algebra} of $\a$, which enjoys a number of universal properties (see \cite[Theorem 10.1.11]{palmer} and \cite[Theorem 10.1.12]{palmer}). For this reason, when every element $a \in \a$ has a bounded universal norm, i.e. $\a_u = \a$, it is said that $\a$ \emph{has an enveloping $C^*$-algebra}.

In general, a $^*$-algebra does not necessarily have an enveloping $C^*$-algebra. Perhaps the most basic example is that of a polynomial $^*$-algebra in a single self-adjoint variable.

We now look at the relation between automatically bounded elements and elements with a bounded universal norm. It is known that every $BG^*$-algebra has an enveloping $C^*$-algebra (\cite[Proposition 10.1.19]{palmer}), and the same proof yields this slightly more general result, that an automatically bounded element has a bounded universal norm:\\

\begin{prop}
\label{autom bounded implies bounded univ norm}
Let $\a$ be a $^*$-algebra. We have that $\a_b \subseteq \a_u$. In particular if $\a$ is a $BG^*$-algebra, then $\a$ has an enveloping $C^*$-algebra.\\
\end{prop}

{\bf \emph{Proof:}} Suppose $a \notin \a_u$. Then there is a sequence of representations $\{\pi_i\}_{i \in \mathbb{N}}$, on Hilbert spaces $\{\mathscr{H}_i\}_{i \in \mathbb{N}}$, such that $\|\pi_i(a)\| \to \infty$. Consider now the inner product space $\mathscr{V}$ defined as the algebraic direct sum
\begin{align*}
\mathscr{V}:= \bigoplus_{i \in \mathbb{N}} \mathscr{H}_i\,,
\end{align*}
and the pre-$^*$-representation $\pi := \bigoplus_{i \in \mathbb{N}}\, \pi_i$ of $\a$ on $L(\mathscr{V})$. It is clear by construction that $\pi(a) \notin B(\mathscr{V})$. Hence, $a \notin \a_b$. \qed\\

\subsection{Preliminaries on Hecke Algebras}

We will mostly follow \cite{krieg} and \cite{schl} in what regards Hecke pairs and Hecke algebras and refer to these references for more details.\\

\begin{df}
 Let $G$ be a group and $\gm$ a subgroup. The pair $(G, \gm)$ is called a \emph{Hecke pair} if every double coset $\gm g\gm$ is the union of finitely many right (and left) cosets. In this case, $\gm$ is also called a \emph{Hecke subgroup} of $G$.\\
\end{df}

Given a Hecke pair $(G, \gm)$ we will denote by $L$ and $R$, respectively, the left and right coset counting functions, i.e.
\begin{align*}
 L(g):= |\gm g \gm / \gm| = [\gm : \gm \cap g \gm g^{-1}] < \infty\\
 R(g) :=|\gm \backslash \gm g \gm| = [\gm : \gm \cap g^{-1} \gm g] < \infty\,.
\end{align*}
We recall that $L$ and $R$ are $\gm$-biinvariant functions which satisfy $L(g) = R(g^{-1})$ for all $g \in G$. Moreover, the function $\Delta: G \to \mathbb{Q^+}$ given by
\begin{align*}
 \Delta(g) := \frac{L(g)}{R(g)}\,,
\end{align*}
is a group homomorphism, usually called the \emph{modular function} of $(G, \gm)$.\\

\begin{df}
 Given a Hecke pair $(G, \gm)$, the \emph{Hecke algebra} $\h(G, \gm)$ is the $^*$-algebra of finitely supported $\mathbb{C}$-valued functions on the double coset space $\gm \backslash G / \gm$ with the product and involution defined by
\begin{align*}
 (f_1*f_2)(\gm g \gm) & := \sum_{h\gm \in G / \gm} f_1(\gm h \gm)f_2(\gm h^{-1}g\gm)\,,\\
f^*(\gm g\gm) & := \Delta(g^{-1}) \overline{f(\gm g^{-1} \gm)}\,.\\
\end{align*}

\end{df}

\begin{rem}
Some authors, including Krieg \cite{krieg}, do not include the factor $\Delta$ in the involution. Here we adopt the convention of Kaliszewski, Landstad and Quigg \cite{schl} in doing so, as it gives rise to a more natural $L^1$-norm. We note, nevertheless, that there is no loss (or gain) in doing so, because these two different involutions give $^*$-isomorphic Hecke algebras. In particular, the question of existence of an enveloping $C^*$-algebra is not perturbed by this.\\
\end{rem}

The Hecke algebra has a natural basis, as a vector space, given by the characteristic functions of double cosets. We will henceforward identify a characteristic function of a double coset $\chi_{\gm g \gm}$ with the double coset $\gm g \gm$ itself. It will be useful to know how to write a product $\gm g \gm * \gm h \gm$ of two double cosets in the unique linear combination of double cosets:\\

\begin{lemma}
\label{unique expression for prod}
 The expression for the product $\gm g \gm * \gm h \gm$ of two double cosets in the unique linear combination of double cosets is given by:
\begin{align*}
 \gm g \gm * \gm h \gm = \sum_{\gm s \gm \in \gm \backslash G / \gm} \frac{L(g)\,C_{g, h}(s)}{L(s)} \,\gm s \gm\,,
\end{align*}
where $C_{g,h}(s):= \# \{w\gm \subseteq \gm h \gm : \gm g w \gm = \gm s \gm\}$.\\
\end{lemma}

{\bf \emph{Proof:}} Let us first check that $C_{g,h}(s)$ is well-defined. It is clear that $C_{g,h}(s)$ does not depend on the representatives $h$ and $s$ of the chosen double cosets, so it remains to verify that it is also independent on $g$. Given any other representative $\beta g \gamma$ of the double coset $\gm g \gm$, with $\beta, \gamma \in \gm$, it is not difficult to see that the map
\begin{align*}
 w\gm \mapsto \gamma^{-1} w\gm\,,
\end{align*}
gives a bijective correspondence between the sets $\{w\gm \subseteq \gm h \gm : \gm g w \gm = \gm s \gm\}$ and $\{u\gm \subseteq \gm h \gm : \gm \beta g \gamma  u \gm = \gm s \gm\}$. Hence we have $C_{\beta g \gamma, h}(s) = C_{g, h} (s)$.

Now, to check the product formula we recall (for example from \cite{schl}) that
\begin{align}
\label{usual prod of double cosets}
 \gm g \gm * \gm h \gm = \sum_{w \gm \in \gm h \gm / \gm} \frac{L(g)}{L(gw)} \,\gm gw \gm\,,
\end{align}
where the sum runs over a set of representatives for left cosets in $\gm h \gm$.
Let us fix a representative $g$ for the double coset $\gm g \gm$ and let $S$ be the set of double cosets $S:= \{\gm g w \gm: w\gm \in \gm h \gm / \gm \}$, i.e. the set of double cosets that appear as summands in (\ref{usual prod of double cosets}). The number of times an element $\gm s \gm \in S$ appears repeated in the sum (\ref{usual prod of double cosets}) is precisely the number $C_{g,h}(s)$. Hence we can write
\begin{align*}
 \gm g \gm * \gm h \gm = \sum_{\gm s \gm \in S} \frac{L(g)\,C_{g,h}(s)}{L(s)} \,\gm s \gm\,.
\end{align*}
Also, if a double coset $\gm r \gm$ does not belong to $S$ we have $C_{g,h}(r) = 0$, thus we get
\begin{align*}
 \gm g \gm * \gm h \gm = \sum_{\gm s \gm \in \gm \backslash G / \gm} \frac{L(g)\,C_{g, h}(s)}{L(s)} \,\gm s \gm\,.
\end{align*}
\qed\\

The reader can find alternative ways of describing the coefficients of this unique linear combination in \cite[Lemma 4.4]{krieg}. In particular, the characterization $(iii)$ of the cited lemma is very similar to the one we just described.\\

\begin{rem}
\label{g gamma h in prod of double cosets}
 A direct computation or Lemma \ref{unique expression for prod} imply that the double cosets that appear in the expression for $\gm g \gm * \gm h \gm$ as a unique linear combination of double cosets are all of the form $\gm g \gamma h\gm$, for some $\gamma \in \gm$. Conversely, all double cosets of the form $\gm g \gamma h\gm$, with $\gamma \in \gm$, appear in this linear combination, because $C_{g,h}(g\gamma h) \neq 0$.\\
\end{rem}

Another basic property of Hecke algebras which we will need is the following: given a Hecke pair $(G, \gm)$ and a subgroup $K$ such that $\gm \subseteq K \subseteq G$, then $(K, \gm)$ is a Hecke pair and $\h(K, \gm)$ is naturally seen as a $^*$-subalgebra of $\h(G, \gm)$. This is a particular case of \cite[Lemma 4.9]{krieg}.\\

\begin{df}
 The \emph{$L^1$-norm} on $\h(G, \gm)$, denoted $\| \cdot \|_{L^1}$, is given by
\begin{align*}
 \| f \|_{L^1} = \sum_{\gm g \gm \in \gm \backslash G / \gm} |f(\gm g \gm)|\, L(g) = \sum_{g \gm \in  G / \gm} |f(\gm g \gm)|
\end{align*}
The completion $L^1(G, \gm)$ of $\h(G, \gm)$ under this norm is a Banach $^*$-algebra.\\
\end{df}

As observed in \cite{bost connes}, \cite{tzanev} and \cite{schl} there are several canonical $C^*$-completions of $\h(G, \gm)$. These are:
 
\begin{itemize}
 \item $C^*_r(G, \gm)$ - Usually called the \emph{reduced Hecke $C^*$-algebra}, is the completion of $\h(G, \gm)$ under the $C^*$-norm arising from a regular representation (see, for example, \cite{tzanev} for details).
 \item $pC^*(\overline{G})p$ - The corner of the full group $C^*$-algebra of the Schlichting completion $(\overline{G}, \overline{\gm})$ of the pair $(G, \gm)$, for the projection $p := \chi_{\overline{\gm}}$. We will not describe this construction here since it is well documented in the literature (\cite{tzanev}, \cite{glockner}, \cite{schl}) and because we will not make use of Schlichting completions in this work.
\item $C^*(L^1(G, \gm))$ - The enveloping $C^*$-algebra of $L^1(G, \gm)$.
\item $C^*(G, \gm)$ - The enveloping $C^*$-algebra (if it exists!) of $\h(G, \gm)$. When it exists, it is usually called the \emph{full Hecke $C^*$-algebra}.
\end{itemize}

These different $C^*$-completions of $\h(G, \gm)$ are related in the following way, through canonical surjective maps:
\begin{align*}
 C^*(G,\gm)  \dashrightarrow C^*(L^1(G, \gm)) \longrightarrow pC^*(\overline{G})p \longrightarrow C^*_r(G, \gm)\,.
\end{align*}

As was pointed out by Hall \cite[Proposition 2.21]{hall}, a Hecke algebra does not need to have an enveloping $C^*$-algebra in general
, with the Hecke algebra of the pair $(SL_2(\mathbb{Q}_p), SL_2(\mathbb{Z}_p))$ being one such example, where $p$ is a prime number and $\mathbb{Q}_p$, $\mathbb{Z}_p$ denote respectively the field of $p$-adic numbers and the ring of $p$-adic integers.

\subsection{Preliminaries on Directed Graphs}

Recall that a \emph{simple directed graph} $\mathcal{G}:=(\b, E)$ consists of a set $\b$, whose elements are called \emph{vertices}, and a subset $E \subseteq \b^2$, whose elements are called \emph{edges}. An edge is thus a pair of vertices $(a,b)$, which we see as directed from $a$ to $b$. Since we are only interested in directed graphs that are simple, i.e. such that there is at most one edge directed from one vertex to another, we will henceforward drop the word \emph{simple} and simply write \emph{directed graph}.

Let us now set some notation. Let $\mathcal{G} := (\b, E)$ be a directed graph. If the ordered pair $(a,b)$ belongs to $E$ we say that $b$ is a \emph{successor} of $a$.\\

\begin{df}
Let $\mathcal{G} := (\b, E)$ be a directed graph. A set of vertices $\mathcal{S} \subseteq \b$ is said to be \emph{co-hereditary} if it contains the successors of all of its elements, i.e. if $a \in S$ and $b \in \b$ is a successor of $a$, then $b \in S$. \\
\end{df}

It is easy to see that an arbitrary intersection of co-hereditary sets is still a co-hereditary set. Hence, we can talk about the co-hereditary set generated by a subset $X \subseteq \b$ of vertices:\\

\begin{df}
Let $\mathcal{G} := (\b, E)$ be a directed graph and $X \subseteq \b$ a set of vertices. The \emph{co-hereditary set generated by} $X$ is the smallest co-hereditary set that contains $X$.\\
\end{df}

Given a directed graph $\mathcal{G} := (\b, E)$ and a set of vertices $X \subseteq \b$, we will denote by $S(X)$ the set of all the successors of all elements of $X$, i.e.
\begin{align*}
 S(X):=\{a \in \b: a\;\; \text{is a successor of}\;\; x,\;\; \text{for some}\;\; x \in X \}\,.
\end{align*}
Similarly, we define the \emph{$n$-th successor set of} $X$ inductively as follows:
\begin{align*}
S^0(X):= X\,, \qquad\qquad S^n(X):=S(S^{n-1}(X))\,, \;\;\text{for}\;\;n \geq 1\,.
\end{align*}

In this way, the $0$-th successor set is simply the set $X$,  the $1$-st successor set is $S(X)$, the $2$-nd successor set is $S(S(X))$, etc. We will often consider $X$ to be a singleton set $X = \{b\}$, and in this case we will use the notation $S(b)$ instead of $S(\{b\})$. The following result follows easily from the definitions:\\

\begin{lemma}
\label{co-her = Sn}
 Let $\mathcal{G} := (\b, E)$ be a directed graph and $X \subseteq \b$ a set of vertices. The co-hereditary set generated by $X$ is the set $\bigcup_{n \in \mathbb{N}_0} S^n(X)$.\\
\end{lemma}

\begin{rem}
 The sets of vertices we are going to consider in our applications will be sets with specific additional structure (for instance, the set of vertices will typically be a basis of a vector space), and we are interested in proving results of the type: all elements of the co-hereditary set generated by $X$ have a certain property $P$. To do so, we use a certain form of ``induction''. Namely, if we prove that all elements of $X$ have the property $P$, and if we prove that the property $P$ is preserved upon taking successors, then by Lemma \ref{co-her = Sn} and the usual induction on $\mathbb{N}$, all elements of the co-hereditary set generated by $X$ will also satisfy $P$.\\
\label{induc}
\end{rem}

\section{Graph Associated with a $^*$-Algebra}
\label{graph assoc with star alg}

Let $A$ be a $^*$-algebra. Suppose that we are given a finite set of elements $\{b_1, \dots, b_n\} \subseteq \a$ satisfying a set of relations of the form
\begin{align}
\label{relations of bi}
 \begin{cases}
  b_1^*b_1 & = \;\;  \lambda_{11}b_1 + \dots + \lambda_{1n}b_n\\
 & \; \vdots  \\ 
 b_n^*b_n & = \;\; \lambda_{n1}b_1 + \dots + \lambda_{nn}b_n\,,
 \end{cases}
\end{align}
where $\lambda_{ij} \in \mathbb{C}$ for each $i,j \in \{1, \dots, n\}$. We claim that the elements $b_1, \dots, b_n$ are automatically bounded, and this fact will pave the way for our study of existence of enveloping $C^*$-algebras:\\

\begin{thm}
\label{theorem for bs}
 Let $A$ be a $^*$-algebra and $\{b_1, \dots, b_n\}\subseteq \a$ a finite set of elements satisfying relations as in (\ref{relations of bi}). Then the elements $b_1, \dots, b_n$ are automatically bounded. In particular they have a bounded universal norm.\\
\end{thm}

In order to prove Theorem \ref{theorem for bs} we will need the following lemma:\\

\begin{lemma}
Let $n \in \mathbb{N}$ and $k_{ij} \in \mathbb{R}_0^+$ for every $1 \leq i,j \leq n$. The set
\begin{align*}
 B:=\{(x_1, \dots, x_n) \in (\mathbb{R}_0^+)^n : x_i^2 \leq k_{i1}x_1+ \dots + k_{in}x_n\qquad \forall 1 \leq i \leq n\}
\end{align*}
is bounded in $\mathbb{R}^n$.\\
\label{lema 2}
\end{lemma}

{\bf \emph{Proof:}} Let us denote by $\beta$ the real number
\begin{align*}
\beta := \sum_{i = 1}^n \sqrt{\sum_{j = 1}^n k_{ij} }\,,
\end{align*}
and let $\widetilde{B}$ be the set defined by
\begin{align*}
\widetilde{B}:= \Big\{(x_1, \dots, x_n) \in (\mathbb{R}_0^+)^n : x_1 + \dots + x_n \leq \beta\, \sqrt{  x_1 + \dots + x_n } \;\Big\}\,,
\end{align*}
We claim that $B \subseteq \widetilde{B}$. To see this, let $(x_1, \dots, x_n) \in B$. We have
\begin{eqnarray*}
 x_1 + \dots + x_n & \leq & \sum_{i=1}^n \sqrt{k_{i1}x_1+ \dots + k_{in}x_n}\\
& \leq & \sum_{i=1}^n \sqrt{\Big( \sum_{j = 1}^n k_{ij} \Big)x_1+ \dots + \Big( \sum_{j = 1}^n k_{ij} \Big)x_n}\\
& = & \sum_{i=1}^n \sqrt{\Big( \sum_{j = 1}^n k_{ij} \Big) \Big(x_1 + \dots + x_n}\Big)\\
& = & \beta\, \sqrt{x_1 + \dots + x_n}\,,
\end{eqnarray*}
and therefore $(x_1, \dots, x_n) \in \widetilde{B}$.

Hence, it is enough to prove that the set $\widetilde{B}$ is bounded. As it is well known,  linear functions in $\mathbb{R}$ grow faster than square roots, thus it is clear that the set
\begin{align*}
Y:= \{x \in \mathbb{R}_0^+ : x \leq \beta \sqrt{x}\} 
\end{align*}
is bounded in $\mathbb{R}$. Let $S:(\mathbb{R}_0^+)^n \to \mathbb{R}$ be the function $S(x_1, \dots, x_n):= \sum_{i=1}^n x_i$. We have that $\widetilde{B} \subseteq S^{-1}(Y)$. Since $S$ is only defined for elements in $(\mathbb{R}_0^+)^n$, the pre-image by $S$ of a bounded set in $\mathbb{R}$ is also a bounded set in $(\mathbb{R}_0^+)^n$. We conclude that $\widetilde{B}$, and therefore $B$, is bounded. \qed\\

{\bf \emph{Proof of Theorem \ref{theorem for bs}:}} Let $\{b_1, \dots, b_n\} \subseteq \a$ be a finite set in $\a$ satisfying relations as in (\ref{relations of bi}) and $B \subseteq (\mathbb{R}^+_0)^n$ the set defined by
\begin{align*}
B:=\{(x_1, \dots, x_n) \in (\mathbb{R}_0^+)^n : x_i^2 \leq |\lambda_{i1}|x_1+ \dots + |\lambda_{in}|x_n\qquad \forall 1 \leq i \leq n\}\,.
\end{align*}
Let $\pi:\a \to L(\mathscr{V})$ be a pre-$^*$-representation and $\xi \in V$ a vector such that $\| \xi \| = 1$. We have that
\begin{eqnarray*}
\| \pi(b_i) \xi \|^2 & = & \langle \pi(b_i^*b_i)\xi, \xi \rangle\;\; \leq \;\; \|\pi(b_i^*b_i)\xi\|\|\xi\|\\
& = & \|\pi(b_i^*b_i)\xi\|\;\; = \;\; \|\sum_{j=1}^n \lambda_{ij}\, \pi(b_j) \xi\|\\
& \leq & \sum_{j=1}^n |\lambda_{ij}|\, \|\pi(b_j)\xi\|\,.
\end{eqnarray*}
Hence it follows that $\big(\|\pi(b_1) \xi\|, \dots, \|\pi(b_n)\xi\|\big) \in B$. Since the definition of the set $B$ is independent of $\pi$ and $\xi$, and since by Lemma \ref{lema 2} we know that $B$ is bounded in $\mathbb{R}^n$, it follows that
\begin{align*}
\sup_{\|\xi \| = 1} \|\pi(b_i) \xi\| < \infty\,,
\end{align*}
i.e. $\pi(b_i) \in B(\mathscr{V})$, for all $1 \leq i \leq n$, and all pre-$^*$-representations $\pi$. Thus, the elements $b_1, \dots, b_n$ are all automatically bounded and therefore have bounded universal norms by Proposition \ref{autom bounded implies bounded univ norm} \qed\\

In practice though, Theorem  \ref{theorem for bs} can be difficult to apply, as in general one is not given a set of elements $\{b_1, \dots, b_n\}$ satisfying the prescribed relations, especially if the structure of the $^*$-algebra $\a$ is not well understood. For this reason we will describe a more algorithmic approach to Theorem \ref{theorem for bs} where the set $\{b_1, \dots, b_n\}$ is not given from the start, but it is instead constructed step-by-step starting from one element $b_1$. This method will be explained through the language of graphs and will be especially useful when applied to Hecke algebras, where knowledge from the Hecke pair can many times be used to show that sets of elements $\{b_1, \dots, b_n\}$ satisfying (\ref{relations of bi}) abound.

Let $\a$ be a $^*$-algebra and $\b$ a basis of $\a$ as a vector space. Given a basis element $b_0 \in \b$ we will denote by $\Phi_{b_0}$ the unique linear functional $\Phi_{b_0}:\a \to \mathbb{C}$ such that

\begin{align}
\Phi_{b_0} (b):=
\begin{cases}
 1, \qquad \text{if}\;\; b = b_0\\
 0, \qquad \text{if}\;\; b \neq b_0
\end{cases}
\end{align}

for every $b \in \b$.\\

\begin{df}
Given a $^*$-algebra $\a$ with a specified basis $\b$, we define its \emph{associated graph} as the directed graph $\mathcal{G}:=(\b, E)$, whose set of vertices is the set $\b$ and whose set of edges is the set
\begin{align*}
E:= \{(a, b) \in \b^2 :\; \Phi_{b} \;(a^* a) \neq 0\}\,.\\
\end{align*}

\end{df}

Thus, given a vertex $a \in \b$, its successors are precisely those basis elements that have non-zero coefficients in the unique expression of $a^*a$ as a linear combination of elements of $\b$, i.e. if
\begin{align*}
a^*a= k_1 b_1 + \dots + k_n b_n\,,
\end{align*}
where each $k_i \in  \mathbb{C}$ is non-zero and the basis elements $b_i$ are all different, then the successors of $a$ are precisely $b_1, \dots , b_n$.\\

\begin{prop}
Let $\a$ be a $^*$-algebra with basis $\b$ and $\mathcal{G}$ its associated graph. If $X \subseteq \b$ is a finite co-hereditary set in $\mathcal{G}$, then all elements of $X$ are automatically bounded.
In particular, all elements of $X$ have a bounded universal norm.\\
\label{claim 1}
\end{prop}

{\bf \emph{Proof:}} Let $X = \{b_1, \dots, b_n\} \subseteq \b$. Since $X$ contains the successors of all its elements, we must necessarily have
\begin{align*}
 b_i^*b_i = \lambda_{i1}b_1 + \dots + \lambda_{in}b_n\,, \qquad\qquad 1 \leq i \leq n \,,
\end{align*}
for some elements $\lambda_{ij} \in \mathbb{C}$ (possibly being zero). It then follows form Theorem \ref{theorem for bs} that all elements $b_1, \dots, b_n$ are automatically bounded and in particular have a bounded universal norm. \qed\\

\begin{cor}
\label{cor principal}
Let $\a$ be a $^*$-algebra, $\b$ a basis for $\a$ and $\mathcal{G}$ its associated graph. If $\a$ is generated as a $^*$-algebra by the elements of the finite co-hereditary sets of $\mathcal{G}$, then $\a$ is a $BG^*$-algebra. In particular $\a$ has an enveloping $C^*$-algebra.\\
\label{corollary 1}
\end{cor}

{\bf \emph{Proof:}} Let $\b_0$ be the set of elements of the finite co-hereditary sets of $\mathcal{G}$. By Proposition \ref{claim 1}, all elements in $\b_0$ are contained in $\a_b$. Since the elements of $\b_0$ generate the $^*$-algebra $\a$, we conclude that $\a = \a_b$, i.e. $\a$ is a $BG^*$-algebra. \qed\\

 We can interpret the above corollary in the following (equivalent) way: suppose we have a $^*$-algebra $\a$ with a basis $\b$. Suppose additionally that we have a particular set $\mathcal{B}_0 \subset \b$ which generates $\a$. If all the elements of $\mathcal{B}_0$ generate finite co-hereditary sets of the associated graph, then $\a$ has an enveloping $C^*$-algebra. Let us now give a couple of immediate examples:\\

\begin{ex}
\label{ex finite dim}
 Let $\a$ be a finite-dimensional $^*$-algebra. If we take any basis $\b$, the associated graph necessarily has finitely many vertices (and edges). Thus, the co-hereditary set generated by any $b \in \b$ is finite.\\
\end{ex}

\begin{ex}
\label{ex discrete group}
 Let $G$ be a discrete group, $\mathbb{C} (G)$ its group algebra with basis $\{\delta_g \in \mathbb{C}(G) : g \in G\}$. Since in the group algebra we have $\delta_g^**\delta_g = \delta_e$, the only successor of $\delta_g$ in the associated graph is $\delta_e$. Since $\delta_e$ is the only successor of itself, the co-hereditary set generated by $\delta_g$ has only two elements, $\delta_g$ and $\delta_e$.\\
\end{ex}

Some non-trivial examples, arising from Hecke algebras, will be computed later in section \ref{classes of Hecke pairs}.\\

\section{A sufficient condition implying the isomorphism $C^*(G, \gm) \cong  C^*(L^1(G, \gm))$}
\label{hecke star alg section}

In Corollary \ref{cor principal} of the previous section we obtained a sufficient condition for a $^*$-algebra to have an enveloping $C^*$-algebra, namely when it is generated by elements of the finite co-hereditary sets (with respect to a given basis). In this section we will improve this result in the case of a Hecke algebra $\h(G, \gm)$: under a suitable assumption we will not only assure an enveloping $C^*$-algebra of $C^*(G, \gm)$ exists, but we will also be able to identify it with $C^*(L^1(G, \gm))$.

Throughout this section and henceforward $(G, \gm)$ will denote a Hecke pair. We will always consider the canonical basis in the Hecke algebra $\h(G, \gm)$, consisting of double cosets $\{\gm g \gm: g \in G\}$. This section is devoted to the proof of the following result:\\

\begin{thm}
\label{enveloping algebra comes from L1}
 Let $(G, \gm)$ be a Hecke pair. If all double cosets generate finite co-hereditary sets, then the enveloping $C^*$-algebra of $\mathcal{H}(G, \gm)$ exists and coincides with $C^*(L^1(G,\gm))$.\\
\end{thm}

In order to give a proof of Theorem \ref{enveloping algebra comes from L1} we will make use of several lemmas.\\

\begin{lemma}
\label{lemma about L1 norm}
 Let $(G, \gm)$ be a Hecke pair and $f_1, f_2 \in \h(G, \gm)$ be two elements such that $f_i(\gm g \gm) \geq 0$ for all $\gm g \gm \in \gm \backslash G / \gm$. The $L^1$-norm satisfies the equality
\begin{align*}
 \| f_1* f_2\|_{L^1} = \|f_1 \|_{L^1}\|f_2\|_{L^1}\,.
\end{align*}
In particular the following equality is also satisfied for any $f \in \h(G, \gm)$ such that $f(\gm g \gm) \geq 0$ for all $\gm g \gm \in \gm \backslash G / \gm$:
\begin{align*}
 \| f^** f\|_{L^1} = \|f \|_{L^1}^2\,.\\
\end{align*}
\end{lemma}

{\bf \emph{Proof:}} We have that
\begin{eqnarray*}
 \|f_1 * f_2 \|_{L^1} & = & \sum_{g\gm \in G / \gm} |(f_1 * f_2)(\gm g \gm)| \\
& = & \sum_{g\gm \in G / \gm} | \sum_{h \gm \in G / \gm} f_1 (\gm h \gm) f_2(\gm h^{-1} g \gm)\;| \\
& = &   \sum_{h \gm \in G / \gm} f_1 (\gm h \gm) \sum_{g\gm \in G / \gm} f_2(\gm h^{-1} g \gm)\\
& = & \sum_{h \gm \in G / \gm} f_1 (\gm h \gm) \sum_{g\gm \in G / \gm} f_2(\gm g \gm)\\
& = & \| f_1 \|_{L^1} \|f_2 \|_{L^1}\,.
\end{eqnarray*}
The second claim in this lemma follows directly from the first statement:
\begin{align*}
 \|f^**f\|_{L^1} = \|f^*\|_{L^1} \|f\|_{L^1} = \|f\|_{L^1}^2\,.
\end{align*}
\qed\\

\begin{lemma}
\label{nonsingularity of matrix lemma}
 Let $n \in \mathbb{N}$ and $A = [a_{ij}]$ be an $n \times n$ matrix, whose entries satisfy: $a_{ii} \in \mathbb{R}^+$ and $a_{ij} \in \mathbb{R}^-_0$ for all $i \neq j$. If there are vectors $\mathbf{d} = (d_1, \dots, d_n)$ and $\mathbf{z} = (z_1, \dots, z_n)$ both in $(\mathbb{R}^+)^n$ satisfying the system
\begin{align}
\label{system}
 A\mathbf{z} = \mathbf{d}\,,
\end{align}
 then $A$ is non-singular.\\
\end{lemma}

{\bf \emph{Proof:}} Let $\mathbf{z} \in (\mathbb{R}^+)^n$ be a solution to the above system. Suppose that $\mathrm{Ker}\, A \neq \{0\}$. Then, the set of solutions to the system (\ref{system}) contains a line $L$. Consider now the set $S$ of all the (finitely many) points which are the intersections of $L$ with the canonical hyperplanes of the form $x_i = 0$, and take a point $\mathbf{y} \in S$ (not necessarily unique) which is closest to $\mathbf{z}$. The point $\mathbf{y}$ is the intersection of $L$ with one of the hyperplanes $x_i=0$, say $x_{i_0} = 0$ with $1 \leq i_0 \leq n$. Since $\mathbf{y} = (y_1, \dots, y_n)$ is in $L$, it is also a solution of the system (\ref{system}) and therefore must satisfy
\begin{align*}
  \sum_{\substack{k = 1 \\ k \neq i_0}}^n a_{i_0 k} y_k = d_{i_0}\,,
\end{align*}
implying that there exists at least one number $y_k$ which is negative. But on the other hand, the open segment between $\mathbf{z}$ and $\mathbf{y}$ lies inside $(\mathbb{R}^+)^n$ because $\mathbf{z} \in (\mathbb{R}^+)^n$ and this segment does not intersect any hyperplane $x_i = 0$ (by choice of the point $\mathbf{y}$). Thus the entries of $\mathbf{y} = (y_1, \dots, y_n)$ are all non-negative, which is a contradiction. Therefore $\mathrm{Ker}\, A = \{0\}$. \qed\\

In preparation for the next lemma we set some notation. Given two vectors $\mathbf{a} =(a_1, \dots, a_n)$ and $\mathbf{b}=(b_1, \dots, b_n) \in \mathbb{R}^n$, we will write $\mathbf{a} \leq \mathbf{b}$ whenever $a_i \leq b_i$ for every $1 \leq i \leq n$. We will denote the zero vector by $\mathbf{0} = (0, \dots, 0)$. Also, given a set of vectors $S \subseteq \mathbb{R}^n$, we will denote by $\mathcal{C}(S)$ the \emph{cone generated by} $S$, i.e. the set of all linear combinations with coefficients in $\mathbb{R}_0^+$ of the elements of $S$.\\

\begin{lemma}
\label{intersec is the point with big coordinates}
Let $n \in \mathbb{N}$ and $A = [a_{ij}]$ be an $n \times n$ matrix whose entries satisfy: $a_{ii} \in \mathbb{R}^+$ and $a_{ij} \in \mathbb{R}^-_0$ for all $i \neq j$. Assume that there are vectors $\mathbf{d} = (d_1, \dots, d_n) > \mathbf{0}$ and $\mathbf{z} = (z_1, \dots, z_n) >\mathbf{0}$ satisfying the system $A\mathbf{z} = \mathbf{d}$. Then, if
\begin{align*}
 A \mathbf{y} \geq \mathbf{0}\,,
\end{align*}
for some $\mathbf{y} \in \mathbb{R}^n$, we must have $\mathbf{y} \geq  \mathbf{0}$.\\
\end{lemma}

{\bf \emph{Proof:}} As we are in the conditions of Lemma \ref{nonsingularity of matrix lemma}, the matrix $A$ is non-singular. First we claim that $\{\mathbf{y}: A\mathbf{y} \geq \mathbf{0}\} = \mathcal{C}(A^{-1}\mathbf{e}_1, \dots, A^{-1}\mathbf{e}_n)$, where $\mathbf{e}_1, \dots, \mathbf{e}_n \in \mathbb{R}^n$ are the canonical unit vectors. The inclusion $\supseteq$ is obvious, while the inclusion $\subseteq$ follows from the fact that if $A \mathbf{y} \geq \mathbf{0}$ then we can write $A\mathbf{y}$ as a positive linear combination of $\mathbf{e}_1, \dots, \mathbf{e}_n$.
 Thus, to prove this lemma it suffices to prove that $A^{-1}\mathbf{e}_k \geq \mathbf{0}$ for every $1 \leq k \leq n$, and we will show this by induction on $n$. The case $n = 1$ is obvious since $a_{11} \in \mathbb{R}^+$. Let us now assume that the result holds for $n-1$, and prove it for $n$. Let $B_k$ be the matrix obtained from $A$ by deleting the $k$-th row and column. Since $A\mathbf{z} = \mathbf{d}$, it follows readily that
\begin{align}
\label{eq with Bk}
 B_k\, \begin{pmatrix}
        z_1\\
        \vdots\\
        z_{k-1}\\
        z_{k+1}\\
        \vdots\\
        z_n
       \end{pmatrix}
 = \begin{pmatrix}
                                 d_1 - a_{1k}z_k\\
                                 \vdots\\
                                 d_{k-1} - a_{k-1\, k}z_k\\
                                 d_{k+1} - a_{k+1\, k}z_k\\
                                 \vdots\\
                                 d_{n} - a_{n\, k} z_k
                                \end{pmatrix}
\end{align}
Since the right hand side of (\ref{eq with Bk}) is a vector in $(\mathbb{R^+})^{n-1}$, and moreover the entries of the matrix $B_k$ satisfy the conditions in the statement of the lemma, we can use the induction hypothesis on the matrix $B_k$. Let $\mathbf{v}:=(v_1, \dots,v_{k-1},v_{k+1},\dots, v_n) \in \mathbb{R}^{n-1}$ be a solution to the equation
\begin{align*}
 B_k \, \mathbf{v} = \begin{pmatrix}
                      d_1\\
                      \vdots\\
                      d_{k-1}\\
                      d_{k+1}\\
                      \vdots\\
                      d_n
                     \end{pmatrix}\,,
\end{align*}
which exists by Lemma \ref{nonsingularity of matrix lemma} (the reason for the chosen indexing of the entries of $\mathbf{v}$ will become clear in the remaining part of the proof). The induction hypothesis tells us that $\mathbf{v} \geq \mathbf{0}$. We also have that
\begin{align*}
 B_k\, \begin{pmatrix}
        z_1\\
        \vdots\\
        z_{k-1}\\
        z_{k+1}\\
        \cdots\\
        z_n
       \end{pmatrix} - B_k\, \mathbf{v} = \begin{pmatrix}
                                 d_1 - a_{1k}z_k\\
                                 \vdots\\
                                 d_{k-1} - a_{k-1\, k}z_k\\
                                 d_{k+1} - a_{k+1\, k}z_k\\
                                 \vdots\\
                                 d_{n} - a_{n\, k} z_k
                                \end{pmatrix} - \begin{pmatrix}
                      d_1\\
                      \vdots\\
                      d_{k-1}\\
                      d_{k+1}\\
                      \vdots\\
                      d_n
                     \end{pmatrix} = \begin{pmatrix}
                      - a_{1k}z_k\\
                      \vdots\\
                       - a_{k-1\, k}z_k\\
                      - a_{k+1\, k}z_k\\
                      \vdots\\
                     - a_{n\, k} z_k
                     \end{pmatrix} \geq \mathbf{0}\,.
\end{align*}
By the induction hypothesis again, we have $z_i - v_i \geq 0$, for $i \neq k$.

 Consider now the vector $\widetilde{\mathbf{v}} \in \mathbb{R}^n$ given by $\widetilde{\mathbf{v}} := (v_1, \dots, v_{k-1},0, v_{k+1}, \dots, v_n)$. We have that
\begin{align*}
 A\, (\mathbf{z} - \widetilde{\mathbf{v}}) = \begin{pmatrix}
                                              d_1\\
                                              \vdots\\
                                              d_{k-1}\\
                                              d_k\\
                                              d_{k+1}\\
                                              \vdots\\
                                              d_n
                                             \end{pmatrix} - \begin{pmatrix}
                                              d_1\\
                                              \vdots\\
                                              d_{k-1}\\
                                              \sum_{i \neq k}^n a_{ki}v_i\\
                                              d_{k+1}\\
                                              \vdots\\
                                              d_n
                                             \end{pmatrix} = \begin{pmatrix}
                                              0\\
                                              \vdots\\
                                              0\\
                                              d_k -\sum_{i \neq k}^n a_{ki}v_i\\
                                              0\\
                                              \vdots\\
                                              0
                                             \end{pmatrix}\,,
\end{align*}
or in other words,
\begin{align*}
 \mathbf{z} - \widetilde{\mathbf{v}} = A^{-1}\, \begin{pmatrix}
                                              0\\
                                              \vdots\\
                                              0\\
                                              d_k -\sum_{i \neq k}^n a_{ki}v_i\\
                                              0\\
                                              \vdots\\
                                              0
                                             \end{pmatrix} = (d_k -\sum_{i \neq k}^n a_{ki}v_i)\; A^{-1}\mathbf{e}_k\,.
\end{align*}

We now notice that $d_k -\sum_{i \neq k}^n a_{ki}v_i > 0$, because all the $a_{ki} \in \mathbb{R}^-_0$ for $k \neq i$, $v_i \geq 0$ as we saw before, and $d_k > 0$. We have already proven that $z_i - v_i \geq 0$, for $i \neq k$, from which it readily follows that $\mathbf{z} - \widetilde{\mathbf{v}} \geq \mathbf{0}$. We can now conclude that
\begin{align*}
 A^{-1} \mathbf{e}_k = \frac{1}{d_k -\sum_{i \neq k}^n a_{ki}v_i}(\mathbf{z} - \widetilde{\mathbf{v}}) \geq \mathbf{0}\,.
\end{align*}
\qed\\

{\bf \emph{Proof of Theorem \ref{enveloping algebra comes from L1}:}} We already know that if all double cosets generate finite co-hereditary sets, then $\mathcal{H}(G,\gm)$ has an enveloping $C^*$-algebra. Thus, it remains to see that this enveloping $C^*$-algebra is the enveloping $C^*$-algebra of $L^1(G,\gm)$, and for this we only need to show that
\begin{align}
\label{inequality universal norm L1 norm}
 \|a\|_u \leq \|a\|_{L^1},
\end{align}
for any $a \in \mathcal{H}(G,\gm)$. Actually we only need to prove (\ref{inequality universal norm L1 norm}) when $a$ is a double coset $a = \gm s\gm$, since the result for a general $a \in \mathcal{H}(G,\gm)$ follows from the following argument: if we write $a$ in the unique linear combination of double cosets, $a = \sum_{i = 1}^n \lambda_i \gm s_i\gm$, then we have
\begin{eqnarray*}
 \|a\|_u & = & \| \sum_{i = 1}^n \lambda_i \gm s_i\gm\|_u \;\; \leq \;\;  \sum_{i=1}^n |\lambda_i| \|\gm s_i\gm\|_u\\
& \leq & \sum_{i = 1}^n |\lambda_i| \| \gm s_i \gm\|_{L^1} \;\; = \;\; \| \sum_{i = 1}^n \lambda_i \gm s_i\gm\|_{L^1}\\
& = & \| a \|_{L^1}\,.
\end{eqnarray*}
 Let therefore $\gm s\gm$ be a double coset and $\{\gm s_1\gm, \dots, \gm s_n\gm\}$ the finite co-hereditary set it generates. By Lemma \ref{unique expression for prod} we have
\begin{align}
\label{expression for prod s_i -1 s_i}
 (\gm s_i\gm)^**\gm s_i\gm = \sum_{j= 1}^n \lambda_{ij}\, \gm s_j\gm\,,
\end{align}
where the coefficients $\lambda_{ij}$ are given by
\begin{align*}
 \lambda_{ij} \;\;:=\;\; \Delta(s_i)\frac{L(s_i^{-1}) C_{s_i^{-1}, s_i}(s_j)}{L(s_j)} \;\;=\;\; \frac{L(s_i) C_{s_i^{-1}, s_i}(s_j)}{L(s_j)}\,.
\end{align*}

Let $B$ be the set
\begin{align*}
 B:= \{(x_1, \dots, x_n) \in (\mathbb{R}_0^+)^n: x_i^2 \leq \lambda_{i1} x_1 + \dots + \lambda_{in} x_n\,, \;\;\;\forall\; 1 \leq i \leq n\}\,.
\end{align*}
Let us also denote by $C$ the subset of $B$ determined by
\begin{align*}
 C:= \{(x_1, \dots, x_n) \in (\mathbb{R}_0^+)^n: x_i^2 = \lambda_{i1} x_1 + \dots + \lambda_{in} x_n\,,\;\;\;\forall\; 1 \leq i \leq n\}\,.
\end{align*}
It follows immediately from the triangle inequality applied to (\ref{expression for prod s_i -1 s_i}) that the universal norm (in fact, any $C^*$-norm) satisfies $\big(\|\gm s_1 \gm\|_u, \dots, \|\gm s_n\gm\|_u\big) \in B$. Moreover, from Lemma \ref{lemma about L1 norm}, the $L^1$-norm satisfies
\begin{eqnarray*}
\|\gm s_i \gm\|_{L^1}^2  & = & \|(\gm s_i \gm)^**\gm s_i \gm \|_{L^1} \;\; = \;\;\sum_{j= 1}^n \lambda_{ij}\, \|\gm s_j\gm\|_{L^1}\,.
\end{eqnarray*}
Thus, $\big(\|\gm s_1 \gm\|_{L^1}, \dots, \|\gm s_n\gm\|_{L^1}\big) \in C$. For ease of reading we will denote by $\mathbf{z}:=(z_1, \dots, z_n)$ the point $\big(\|\gm s_1\gm\|_{L^1}, \dots, \|\gm s_n\gm\|_{L^1}\big)$. The idea for the remaining part of the proof is to argue that $\mathbf{z} \in C$ is the point with the largest coordinates in the whole set $B$.

For each $1 \leq i \leq n$ let $g_i: (\mathbb{R}_0^+)^n \to \mathbb{R}$ be the function
\begin{align*}
 g_i(x_1, \dots, x_n) := x_i^2 - \sum_{j= 1}^n \lambda_{ij} x_j\,.
\end{align*}
The tangent hyperplane to the graph of $g_i$ at the point $(z_1, \dots, z_n)$ is given by the equation
\begin{align*}
 (2z_i - \lambda_{ii})(x_i - z_i) - \sum_{\substack{j= 1 \\ j \neq i}}^n \lambda_{ij} (x_j - z_j) = 0\,,
\end{align*}
which, using the fact that $(z_1, \dots, z_n)$ is a zero of $g_i$, we can reduce to
\begin{align}
\label{eq of planes}
 (2 z_i - \lambda_{ii})x_i - \sum_{\substack{j= 1 \\ j \neq i}}^n \lambda_{ij}\, x_j = z_i^2\,.
\end{align}
We claim that $2 z_i - \lambda_{ii} > 0$. To see this we notice that
\begin{eqnarray*}
 2z_i - \lambda_{ii} & = & 2\|\gm s_i \gm \|_{L^1} - \frac{L(s_i)C_{s_i^{-1},s_i}(s_i)}{L(s_i)} \;\; = \;\; 2L(s_i) - C_{s_i^{-1},s_i}(s_i)\\
& \geq & 2L(s_i) - L(s_i) \;\; = \;\; L(s_i) \;\; > \;\;0\,.
\end{eqnarray*}

Let us now take $A=[a_{ij}]$ to be the $n \times n$ matrix whose entries are given by $a_{ij} := - \lambda_{ij}$ for $i \neq j$, and $a_{ii} := 2z_i - \lambda_{ii}$, thus $a_{ij} \in \mathbb{R}^-_0$ for $i\neq j$ and $a_{ii} \in \mathbb{R}^+$. We can easily see from (\ref{eq of planes}) that $A\mathbf{z} = \mathbf{z}^2$, where $\mathbf{z}^2 = (z_1^2, \dots, z_n^2)$. Consider now the set $W$ defined by
\begin{align*}
 W:=\big\{\mathbf{x} \in (\mathbb{R}^+_0)^n : A \mathbf{x} \leq \mathbf{z}^2\}\,.
\end{align*}
 We claim that $W$ contains the set $B$. To see this, let $(y_1, \dots, y_n) \in B$. We then have
\begin{eqnarray*}
 -z_i^2 + (2 z_i - \lambda_{ii})y_i - \sum_{\substack{j= 1 \\ j \neq i}}^n \lambda_{ij}\, y_j & = & -z_i^2 +2z_iy_i - \lambda_{ii}y_i - \sum_{\substack{j= 1 \\ j \neq i}}^n \lambda_{ij}\, y_j\\
& = & -(y_i - z_i)^2 + y_i^2 - \sum_{j= 1 }^n \lambda_{ij}\, y_j\\
& \leq &  y_i^2 - \sum_{j= 1 }^n \lambda_{ij}\, y_j\\
& \leq & 0\,,
\end{eqnarray*}
which implies that
\begin{align*}
 (2 z_i - \lambda_{ii})y_i - \sum_{\substack{j= 1 \\ j \neq i}}^n \lambda_{ij}\, y_j \leq z_i^2\,,
\end{align*}
and thus $(y_1, \dots, y_n) \in W$. In other words, if $\mathbf{y} \in B$, then $A\mathbf{y} \leq \mathbf{z}^2$. We can rewrite this inequality as:
\begin{eqnarray*}
 A\mathbf{y} \leq \mathbf{z}^2 & \Leftrightarrow & 0 \leq \mathbf{z}^2 - A\mathbf{y} \;\; \Leftrightarrow \;\; 0 \leq A(\mathbf{z} - \mathbf{y})\,.
\end{eqnarray*}
 Noting that we are under the conditions of Lemma \ref{intersec is the point with big coordinates}, because the entries of $A$ satisfy the required conditions and $A\mathbf{z} = \mathbf{z}^2$, we conclude that $0 \leq \mathbf{z} - \mathbf{y}$, i.e. $\mathbf{y} \leq \mathbf{z}$. Thus, we conclude that $\mathbf{z}$ has bigger coordinates than any other point in $B$.

As we know, we have $(\| \gm s_1\gm \|_u, \dots, \|\gm s_n\gm\|_u) \in B$, so by the above we must have $\| \gm s_i \gm\|_u \leq z_i = \|\gm s_i \gm\|_{L^1}$ for any $1 \leq i \leq n$. Thus, in particular, $\| \gm s\gm \|_u \leq \| \gm s\gm \|_{L^1}$, for the initial double coset $\gm s \gm$. Since all double cosets generate finite co-hereditary sets we conclude that this inequality holds for any double coset $\gm s\gm$, and as we explained in the beginning of the proof, this implies that the enveloping $C^*$-algebra of $\mathcal{H}(G, \gm)$ is $C^*(L^1(G,\gm))$. \qed\\

\section{Methods for Hecke Algebras}
\label{methods section}

The basis of our study of enveloping $C^*$-algebras of Hecke algebras will be Corollary \ref{cor principal} and Theorem \ref{enveloping algebra comes from L1}. Our goal is to apply these results to several classes of Hecke pairs,  but so far we have not given any hint on how to actually ensure that a given double coset generates a finite co-hereditary set. The objective of this section is to provide some tools, based on iterated commutators, to help us accomplish this task.

Given a group $G$ we will denote by $[s,t]$ the commutator of $s, t \in G$, i.e.
\begin{align*}
 [s,t]:= s^{-1}t^{-1}st\,.
\end{align*}
More generally, given elements $s_1, \dots, s_n \in G$ we will denote by $[s_1, \dots, s_n]$ the iterated commutator defined inductively by
\begin{align*}
[s_1, \dots, s_n] := [[s_1, \dots, s_{n-1}], s_n]\,. 
\end{align*}

Let us now return to Hecke pairs $(G, \gm)$. We will be mostly interested in commutators of the form $[g, \gamma_1, \dots, \gamma_n]$, where $g \in G$ and $\gamma_1, \dots, \gamma_n \in \gm$, and the reason for that is given by the following result:\\

\begin{prop}
 Let $(G,\gm)$ be a Hecke pair and $g \in G$. Let $\{\gm x_n \gm\}_{n \in \mathbb{N}_0}$ be a sequence of double cosets satisfying the properties:
\begin{itemize}
 \item[i)] $\gm x_0 \gm = \gm g \gm$,
 \item[ii)] $\gm x_{n+1} \gm$ is a successor of $\gm x_n \gm$, for all $n \geq 0$.
\end{itemize}
 Then, there exists a sequence $\{\gamma_n\}_{n \in  \mathbb{N}} \subseteq \gm$ such that
\begin{align*}
 \gm x_n \gm = \gm [g, \gamma_1, \dots, \gamma_n]\gm\,,
\end{align*}
 for all $n \geq 1$. In particular, all elements in $S^n(\gm g\gm)$ have a representative of the form $\gm [g, \gamma_1, \dots, \gamma_n] \gm$, for some $\gamma_1, \dots, \gamma_n \in \gm$.\\
\label{levels and series}
\end{prop}

{\bf \emph{Proof:}} We will choose such a sequence $\{\gamma_n\}_{n \in \mathbb{N}}$ inductively on $n \in \mathbb{N}$. Suppose $n = 1$. Since $\gm x_1 \gm$ is a successor of $\gm g\gm$, it must be of the form $\gm x_1 \gm = \gm g^{-1}\gamma g\gm$ for some $\gamma \in \gm$ (see Remark \ref{g gamma h in prod of double cosets}). Now we notice that
\begin{align*}
 g^{-1}\gamma g = g^{-1}\gamma g\gamma^{-1}\gamma = [g, \gamma^{-1}]\gamma\,.
\end{align*}
Hence, we have 
\begin{align*}
 \gm x_1 \gm = \gm g^{-1}\gamma g\gm = \gm [g, \gamma^{-1}]\gamma\gm = \gm [g, \gamma^{-1}]\gm\,.
\end{align*}
Choosing $\gamma_1 := \gamma^{-1}$ yields the desired result.

Now let us suppose that there exist elements $\gamma_1, \dots, \gamma_n \in \gm$ such that  $\gm x_k \gm = \gm [g,\gamma_1, \dots, \gamma_k] \gm$, for every $1 \leq k \leq n$. Then, since $\gm x_{n+1}\gm$ is a successor of $\gm x_n \gm = \gm [g, \gamma_1, \dots, \gamma_n] \gm$, we can write
\begin{align*}
 \gm x_{n +1} \gm = \gm [g, \gamma_1, \dots, \gamma_n]^{-1} \gamma [g, \gamma_1, \dots, \gamma_n]\gm\,,
\end{align*}
for some $\gamma \in \gm$ (again by Remark \ref{g gamma h in prod of double cosets}). We have
\begin{eqnarray*}
 \gm x_{n +1} \gm & = & \gm [g, \gamma_1, \dots, \gamma_n]^{-1} \gamma [g, \gamma_1, \dots, \gamma_n] \gamma^{-1}\gm\\
 & = & \gm [g, \gamma_1, \dots, \gamma_n, \gamma^{-1}]\gm\,.
\end{eqnarray*}
Choosing $\gamma_{n+1} := \gamma^{-1}$ yields the desired result for $\gm x_{n+1} \gm$.

Hence, since we can extend any finite sequence $\gamma_1, \dots, \gamma_n$ satisfying the stated conditions to a sequence $\gamma_1, \dots, \gamma_n, \gamma_{n+1}$ still satisfying the stated conditions, it follows that there must be an infinite sequence $\{ \gamma_n \}_{n \in \mathbb{N}}$ with the desired requirements. \qed\\

We will now establish a sufficient condition to ensure the finiteness of the co-hereditary set generated by an element $\gm g\gm$ based on the iterated commutators we considered above:\\

\begin{thm}
 Let $(G,\gm)$ be a Hecke pair and $g \in G$. Suppose that for any sequence of elements $\{\gamma_k\}_{k \in \mathbb{N}} \subseteq \gm$ the total number of double cosets
\begin{align*}
 \#\big\{\gm [g, \gamma_1, \dots, \gamma_n] \gm : n \in \mathbb{N} \big\}
\end{align*}
is finite. Then $\gm g \gm$ generates a finite co-hereditary set.\\
\label{grafo e serie}
\end{thm}

{\bf \emph{Proof:}} Suppose the co-hereditary set generated by $\gm g\gm$ is infinite. Then, there must exist a sequence $\{\gm x_n\gm\}_{n \in \mathbb{N}_0}$ such that
\begin{itemize}
 \item[i)] $\gm x_0 \gm = \gm g \gm$,
 \item[ii)] $\gm x_{n+1} \gm$ is a successor of $\gm x_n \gm$, for all $n \geq 0$.
 \item[iii)] $\gm x_{n+1}\gm \notin \bigcup_{i=0}^n S^i(\gm g\gm)$, for all $n \geq 0$.
\end{itemize}
In particular, we have that $\gm x_i\gm \neq \gm x_j\gm$ for $i \neq j$, implying that the set $\{ \gm x_n \gm: n \in \mathbb{N} \}$ is infinite.

 By Proposition \ref{levels and series} there exists a sequence $\{ \gamma_n\}_{n \in \mathbb{N}} \subseteq \gm$ such that $\gm x_n \gm = \gm [g, \gamma_1, \dots, \gamma_n]\gm$ for all $n \geq 1$. But, by assumption, the number of double cosets in $\{\gm [g, \gamma_1, \dots, \gamma_n] \gm : n \in \mathbb{N}\}$ is finite. Thus we arrive at a contradiction and therefore the co-hereditary set generated by $\gm g \gm$ must be finite. \qed\\

\begin{cor}
 Let $(G,\gm)$ be a Hecke pair and $g \in G$. If one of the following conditions holds, then $\gm g\gm$ generates a finite co-hereditary set:
\begin{itemize}
 \item[a)] For every sequence $\{\gamma_n\}_{n \in \mathbb{N}} \subseteq \gm$ there exists a finite set $F \subseteq G$ and $N_0 \in \mathbb{N}$ such that $[g, \gamma_1, \dots, \gamma_k] \in F$ for all $k \geq N_0$.
 \item[b)] For every sequence $\{\gamma_n\}_{n \in \mathbb{N}} \subseteq \gm$ there exists $N \in \mathbb{N}$ such that $[g, \gamma_1, \dots, \gamma_{N}] \in \gm $.
\end{itemize}

\label{D_k finito}
\end{cor}

{\bf \emph{Proof:}} The result follows directly from Theorem \ref{grafo e serie}. For $a)$ we notice that we can write $ \{ \gm [g, \gamma_1, \dots, \gamma_n]\gm: n \in \mathbb{N} \}$ as the union of the two finite sets $\{ \gm [g, \gamma_1, \dots, \gamma_n]\gm: n < N_0 \}$ and $\{ \gm [g, \gamma_1, \dots, \gamma_n]\gm: n \geq N_0 \}$.

For $b)$, one can easily show, by induction, that $[g, \gamma_1, \dots, \gamma_n] \in \gm$, for any $n \geq N$. Thus, we have
\begin{align*}
 \big\{\gm [g, \gamma_1, \dots, \gamma_n] \gm : n \in \mathbb{N} \big\} = \big\{\gm [g, \gamma_1, \dots, \gamma_n] \gm : n \leq N \big\}\,,
\end{align*}
which is a finite set. \qed\\

There are different classes of Hecke pairs that satisfy conditions $a)$ and  $b)$ of the above corollary. As we shall see in more detail in the next section, condition $a)$ is satisfied by groups satisfying certain generalized nilpotency properties, whereas $b)$ is satisfied when $\gm$ is a subnormal subgroup of $G$, for example.\\

\section{Classes of Hecke Pairs}
\label{classes of Hecke pairs}

We will now use the methods developed in the previous sections to study the existence of enveloping $C^*$-algebras for several classes of Hecke algebras. Many of the well known results about the existence of a full Hecke $C^*$-algebra for some classes of Hecke pairs will be recovered in a unified approach and some new classes will also be described. The isomorphism $C^*(G, \gm) \cong C^*(L^1(G, \gm))$ will also be established in many of the considered classes.

 It should also be noted that all the classes of Hecke algebras considered here are in fact $BG^*$-algebras, since our methods can be traced back to Corollary \ref{cor principal}, but since the focus is mostly on the existence of $C^*(G, \gm)$ we will not mention this in every case.

This section is organized as follows: the classes of Hecke pairs from \ref{finite index class} to \ref{protonormal class} have been studied in the operator algebraic literature and results about the corresponding full Hecke $C^*$-algebras are known. The results about the remaining classes, \ref{subnormal section} to \ref{loc finite classes}, are essentially new, with the results for the classes \ref{subnormal section}, \ref{ascendant section} and \ref{finite conj classes} generalizing known results in the literature.

The classes we consider are presumably all different (in the sense of containment), with the notable exceptions of \ref{subnormal section} which is a particular case of \ref{ascendant section}, and \ref{finite index class} which is a particular case of \ref{finite conj classes}.

We would like to remark that the results discussed in this section illustrate how our methods apply for natural classes of Hecke pairs and that we have not, by any means, exhausted all the possible classes of Hecke pairs one can study through these methods.\\

\subsection{$\gm$ has Finite Index in $G$}
\label{finite index class}

When $\gm$ has finite index in $G$, the pair $(G,\gm)$ is automatically a Hecke pair, and the Hecke $^*$-algebra is finite dimensional (actually, $\h(G, \gm)$ is finite dimensional if and only if $\gm $ has finite index in $G$). As we have seen in Example \ref{ex finite dim}, the co-hereditary set generated by a double coset is finite because the graph of $\h(G, \gm)$ is itself finite. Hence, Theorem \ref{enveloping algebra comes from L1} tells us that $C^*(G, \gm)$ exists and $C^*(G, \gm) \cong C^*(L^1(G, \gm))$.

Of course this example, investigated by Hall \cite[Section 4.2]{hall}, is well-known and completely understood,  because a finite dimensional $^*$-algebra is automatically complete for any $^*$-algebra norm. Hence we necessarily have
\begin{align*}
 C^*(G, \gm) \cong C^*(L^1(G, \gm)) \cong pC^*(\overline{G})p \cong C^*_r(G, \gm)\,,
\end{align*}
and all these $C^*$-algebras are isomorphic to $\h(G, \gm)$, without having to invoke our Theorem \ref{enveloping algebra comes from L1}.

\subsection{$(G,\gm)$ is Directed}

Recall that $(G,\gm)$ is said to be \emph{directed} if $G=T^{-1}T$, where
\begin{align*}
T:=\{t \in G : \gm \subseteq t\gm t^{-1}\}\,.
\end{align*}
Directed Hecke pairs have been widely studied in the literature (\cite{brenken}, \cite{hall}, \cite{larsen raeb}, \cite{laca larsen}, \cite{grp ext}, \cite{schl}, for example), in particular because of their association with the theory of semigroup $C^*$-crossed products. It is known that when $(G,\gm)$ is directed the Hecke algebra has an enveloping $C^*$-algebra and moreover one has
\begin{align*}
 C^*(G, \gm) \cong C^*(L^1(G, \gm)) \cong pC^*(\overline{G})p\,,
\end{align*}
(see, for example,  \cite[Theorem 7.4]{schl}).

 With our methods we can show that $C^*(G, \gm)$ exists, since the Hecke algebra is in fact generated by finite co-hereditary sets. To see this, we first notice that, for $t \in T$, we have $\gm t \gm= t\gm$. Hence, we also have
\begin{align}
(\gm s\gm)^**\gm t\gm= \gm s^{-1}t\gm
\label{equality dir}
\end{align}
for every $s, t \in T$, which means that the Hecke $^*$-algebra is generated by the set of double cosets $\{\gm t\gm: t \in T\}$. Taking $s=t$ in equality (\ref{equality dir}) we see that
\begin{align*}
(\gm t \gm)^**\gm t \gm= \gm
\end{align*}
Thus, the only successor of the double coset $\gm t \gm $ is $\gm$. Since $\gm$ is the only successor of itself, it follows that the co-hereditary set generated by $\gm t \gm$ has only two elements, $\gm t\gm $ and $\gm$, and is therefore finite. We conclude that $\mathcal{H}(G,\gm)$ is generated by finite co-hereditary sets and therefore $C^*(G, \gm)$ exists by Corollary \ref{cor principal}.

\subsection{Iwahori Hecke Algebras}
\label{iwahori hecke class}

Let $(G,\gm)$ be a Hecke pair such that $\mathcal{H}(G,\gm)$ is an Iwahori Hecke algebra (see \cite[Definition 5.12]{hall} for a precise definition of this concept). Sets of generators and relations have been given for this class of Hecke algebras, but for our purposes we will only need to know that:
\begin{enumerate}
 \item There is a set $S \subseteq G$ of elements of order two such that $\mathcal{H}(G,\gm)$ is generated (as a $^*$-algebra) by $\gm$ and the double cosets $\gm s\gm$, with $s \in S$.
 \item for every $s \in S$ the following relation holds:
\begin{align*}
(\gm s\gm)^2 = L(s)\gm + (L(s)-1)\gm s \gm\,.
\end{align*}
\end{enumerate}

For the remaining relations in $\mathcal{H}(G,\gm)$, of which we will not make any use in this work, we refer the reader to Hall's thesis \cite[Section 5.3.1]{hall}.

It was proven by Hall \cite[Proposition 2.24]{hall}, through an estimate on the spectral radius of certain elements, that an Iwahori Hecke algebra has an enveloping $C^*$-algebra (actually Hall proved this for the case $(SL_n(\mathbb{Q}_p), B)$, with $B \subseteq SL_n(\mathbb{Q}_p)$ an Iwahori subgroup, but her proof is completely general).

We can also conclude this from our methods, by proving that $\mathcal{H}(G,\gm)$ is generated by finite co-hereditary sets. By point $1)$ we only need to see that each double coset $\gm s \gm$ with $s \in S$ generates a finite co-hereditary set. So let $\gm s \gm \in \mathcal{H}(G,\gm)$ with $s \in S$. Since $s$ has order two we see that $\gm s \gm$ is self-adjoint and therefore relation $2)$ can be rewritten as
\begin{align*}
(\gm s\gm )^** \gm s\gm = L(s) \gm + (L(s) - 1) \gm s\gm
\end{align*}
Hence, the successors of $\gm s\gm$ are only $\gm$ and $\gm s\gm$ itself. Thus, the co-hereditary set generated by $\gm s\gm$ has only two elements, $\gm$ and $\gm s\gm$, and is therefore finite. We conclude that $\mathcal{H}(G,\gm)$ is generated by finite co-hereditary sets, and is therefore a $BG^*$-algebra and has an enveloping $C^*$-algebra.\\

\begin{rem}
 By a result of Hall \cite[Theorem 6.10]{hall} and a result of Kaliszewski, Landstad and Quigg \cite[Corollary 6.11]{schl} it is known that, for $G = SL_2(\mathbb{Q}_p)$ and $\gm$ an Iwahori subgroup, we necessarily have
\begin{align*}
 C^*(G, \gm) \cong C^*(L^1(G, \gm)) \cong pC^*(\overline{G})p\,.
\end{align*}
 The analogous result for $SL_n(\mathbb{Q}_p)$ with $n \geq 3$ is still open, as far as we know.\\
\end{rem}

\subsection{$\gm$ is a Protonormal Subgroup of $G$}
\label{protonormal class}

We recall that $\gm$ is a \emph{protonormal} subgroup of $G$ (in the sense of Exel \cite{exel}), if for every $s \in G$ we have
\begin{align*}
\gm s^{-1}\gm s = s^{-1}\gm s\gm\,.
\end{align*}
Subgroups with this property are also called \emph{conjugate permutable subgroups} in the literature.

It was proven by Exel (\cite[Proposition 12.1]{exel}) that when $\gm$ is a protonormal subgroup of $G$ the enveloping $C^*$-algebra $C^*(G, \gm)$ exists. Moreover, it is completely clear from his proof that $C^*(G, \gm) \cong C^*(L^1(G, \gm))$, since the bound he uses for the universal norm is actually the $L^1$-norm. Our methods can also recover this result, because in fact any double coset $\gm g \gm$ generates a finite co-hereditary set. We will actually prove that the co-hereditary set generated by $\gm g \gm$ consists only of $\gm g \gm$ and $S(\gm g \gm)$  and is therefore finite. In other words, we will prove that
\begin{align*}
 S^n(\gm g \gm) \subseteq S(\gm g \gm)\,,
\end{align*}
 for every $n \in \mathbb{N}$. It suffices to prove that $S^2(\gm g \gm) \subseteq S(\gm g \gm)$. The elements of $S^2(\gm g \gm)$ are of the form $\gm [g, \gamma_1, \gamma_2] \gm$, where $\gamma_1, \gamma_2 \in \gm$,  by Proposition \ref{levels and series}. We have that
\begin{eqnarray*}
 [g, \gamma_1, \gamma_2] & = & [g, \gamma_1]^{-1} \gamma_2^{-1} [g, \gamma_1]\gamma_2\\
 & = & \gamma_1^{-1}g^{-1}(\gamma_1g\gamma_2^{-1} g^{-1}) \gamma_1^{-1} g \gamma_1 \gamma_2\,.
\end{eqnarray*}
Since $\gm$ is a protonormal subgroup there exist $\theta, \omega \in \gm$ such that $\gamma_1 g \gamma_2^{-1} g^{-1} = g \theta g^{-1} \omega$. Thus, we get
\begin{eqnarray*}
 [g, \gamma_1, \gamma_2] & = & \gamma_1^{-1}g^{-1}(g \theta g^{-1} \omega) \gamma_1^{-1} g \gamma_1 \gamma_2\\
 & = & \gamma_1^{-1}\theta g^{-1} \omega \gamma_1^{-1} g \gamma_1 \gamma_2\,,
\end{eqnarray*}
 and therefore
\begin{eqnarray*}
 \gm [g, \gamma_1, \gamma_2] \gm & = & \gm \gamma_1^{-1}\theta g^{-1} \omega \gamma_1^{-1} g \gamma_1 \gamma_2 \gm\\
 & = & \gm g^{-1} \omega \gamma_1^{-1} g \gm\,.
\end{eqnarray*}
By Remark \ref{g gamma h in prod of double cosets}, $\gm g^{-1} \omega \gamma_1^{-1} g \gm \in S(\gm g \gm)$. This finishes the proof.\\

\subsection{$\gm$ is Subnormal in $G$}
\label{subnormal section}

Hecke pairs $(G,\gm)$ in which $\gm$ is normal in a normal subgroup of $G$ have been widely studied in the literature, in particular when $G$ is a semi-direct product (\cite{brenken}, \cite{larsen raeb}, \cite{laca larsen}, \cite{schl}), and it is known that in this case $\h(G, \gm)$ has an enveloping $C^*$-algebra and moreover
\begin{align*}
 C^*(G, \gm) \cong C^*(L^1(G, \gm)) \cong pC^*(\overline{G})p\,,
\end{align*}
(see, for example, \cite[Theorem 6.13]{schl}).

 We are now going to prove that when $\gm$ is a subnormal subgroup of $G$, $C^*(G, \gm)$ exists and $C^*(G, \gm) \cong C^*(L^1(G, \gm))$. Recall that $\gm$ is \emph{subnormal} in $G$ if there are subgroups $H_0, H_1, \dots, H_n$ such that
\begin{align*}
\gm = H_n \unlhd H_{n-1} \unlhd \dots \unlhd H_0 =G\,,
\end{align*}
where the notation $H_{i+1} \unlhd H_i$ means that $H_{i+1}$ is a normal subgroup of $H_i$.

We claim that when $\gm$ is subnormal in $G$, all double cosets $\gm s \gm$ generate finite co-hereditary sets. To see this we will use Corollary \ref{D_k finito}. Let $s \in G$ and $\{ \gamma_k\}_{k \in \mathbb{N}} \subseteq \gm$. We will prove by induction that $[s, \gamma_1, \dots, \gamma_k] \in H_k$ for $1 \leq k \leq n$. For $k =1$ this follows from the following observation:
\begin{align*}
 [s, \gamma_1]  =  s^{-1}\gamma_1^{-1}s \gamma_1  \in   s^{-1}\gm s \gamma_1  \subseteq  s^{-1}H_1 s \gamma_1  =  H_1 \gamma_1  =  H_1\,.
\end{align*}
Now, let us prove that $k \Rightarrow k +1$. For simplicity, let us write $x_k := [s, \gamma_1, \dots, \gamma_k]$, which by induction hypothesis is an element of $H_k$. Thus, we have
\begin{eqnarray*}
 [s, \gamma_1, \dots, \gamma_k, \gamma_{k+1}] & = & [x_k, \gamma_{k+1}] \;\; \in \;\; x_k^{-1}\gm x_k \gamma_{k+1}\\
& \subseteq & x_k^{-1}H_{k+1} x_k \gamma_{k+1} \;\; = \;\; H_{k+1} \gamma_{k+1}\\
& = & H_{k+1}\,.
\end{eqnarray*}

Thus, for any sequence $\{\gamma_k\}_{k \in \mathbb{N}}$ we have $[s, \gamma_1, \dots, \gamma_n] \in \gm$, which by Corollary \ref{D_k finito} $b)$ implies that $\gm s \gm$ generates a finite co-hereditary set. Since this is true for all double cosets, Theorem \ref{enveloping algebra comes from L1} tells us that $C^*(G, \gm)$ exists and $C^*(G, \gm) \cong C^*(L^1(G, \gm))$.\\

\begin{rem}
 It is known that any subgroup $\gm$ of a nilpotent group $G$ is necessarily a subnormal subgroup (see, for example, \cite[\S 62]{kurosh}). Hence already from this we can conclude that the Hecke algebra of any Hecke pair $(G,\gm)$, with $G$ a nilpotent group, has an enveloping $C^*$-algebra (which coincides with $C^*(L^1(G, \gm))$). In fact, this holds for any group $G$ whose subgroups are all subnormal. Groups with this property form a class that strictly contains the class of nilpotent groups (\cite[Theorem 6.11]{rob}).
We will prove similar results for other classes of groups which strictly generalize the class of nilpotent groups.\\
\label{remark subnormal}
\end{rem}

\begin{ex}
 Let $G$ be the group of $n \times n$ upper triangular matrices with $1$'s on the diagonal and with entries in $\mathbb{Q}$ and let $\gm$ be the subgroup of those matrices with entries in $\mathbb{Z}$. It can be checked, although we will not do so here, that $(G, \gm)$ forms a Hecke pair. The subgroup $\gm$ is subnormal with
\begin{align*}
 \gm = H_n \unlhd H_{n-1} \unlhd \dots \unlhd H_1 =G\,,
\end{align*}
where $H_k$ is the subgroup of matrices in $G$ whose first $k-1$ upper diagonals have entries in $\mathbb{Z}$. The group $G$ is nilpotent and its $3 \times 3$ version is the rational Heisenberg group discussed in \cite[Example 11.7]{schl}.\\
\end{ex}

\subsection{$\gm$ is Ascendant in $G$}
\label{ascendant section}

Recall that $\gm$ is said to be \emph{ascendant} in $G$ if there is a normal series $\{H_i\}_{i \in \mathbb{N}_0}$,
\begin{align*}
 \gm = H_0 \unlhd H_1 \unlhd \dots \unlhd H_i \unlhd \dots 
\end{align*}
that ends in the group $G$, in the sense that $\bigcup_{i \in \mathbb{N}_0} H_i = G$. Of course, the series is finite precisely when $\gm$ is subnormal in $G$.

We will now prove that if $\gm$ is ascendant in $G$, then every double coset generates a finite co-hereditary set, therefore implying that  $C^*(G, \gm)$ exists and is isomorphic to $C^*(L^1(G, \gm))$.

 Let $\gm s\gm $ be any double coset in $\mathcal{H}(G,\gm)$, with representative $s \in G$. Since $\gm$ is ascendant, $s$ must belong to one of the subgroups $H_n$, with $n \in \mathbb{N}_0$. Of course, $\gm$ is a subnormal subgroup of $H_n$, and as we saw in the subnormal case, this implies that the co-hereditary set generated by $\gm s\gm$ is necessarily finite.\\

\subsection{$\gm$ has Finitely Many Conjugates in $G$}
\label{finite conj classes}

Suppose $\gm$ has finitely many conjugates in $G$, or equivalently, the normalizer of $\gm$ has finite index in $G$. Then, $C^*(G, \gm)$ exists and $C^*(G, \gm) \cong C^*(L^1(G, \gm))$ because any double coset generates a finite co-hereditary set. To see this, let $\gm g \gm$ be a double coset and let $g_1^{-1} \gm g_1, \dots, g_n^{-1} \gm g_n$ be the conjugates of $\gm$. With the possible exception of $\gm g \gm$ itself, any element in the co-hereditary set generated by $\gm g \gm$ is a successor of another element. Hence, by Remark \ref{g gamma h in prod of double cosets}, any such element is of the form
\begin{align*}
 \gm x^{-1} \gamma x\gm\,,
\end{align*}
where $x \in G$ and $\gamma \in \gm$. We can then write $ x^{-1} \gamma x =  g_i^{-1} \theta g_i$, for some $i \in\{1, \dots, n\}$ and $\theta \in \gm$, and therefore  $\gm x^{-1} \gamma x \gm = \gm  g_i^{-1} \theta g_i \gm$. Thus, apart possibly from $\gm g\gm$, all elements in the co-hereditary set generated by $\gm g\gm$ are successors of some $\gm g_i \gm$, with $1 \leq i \leq n$, by Remark \ref{g gamma h in prod of double cosets} again. Thus, this co-hereditary set must be finite.\\

\subsection{$G$ is Finite-by-Nilpotent}
\label{finite-by-nil}

Recall that a group $G$ is called \emph{nilpotent} if its lower central series stabilizes at $\{e\}$ after finitely many steps, i.e. if the normal series defined inductively by
\begin{align*}
 G_0:= G\,, \qquad\qquad G_{n+1}:= [G_n, G]\,,
\end{align*}
is such that $G_k=\{e\}$, for some $k \in \mathbb{N}$.

Recall also that a group $G$ is said to be \emph{finite-by-nilpotent} if $G$ has a finite normal subgroup $K$ such that $G/K$ is nilpotent, i.e. if $G$ is an extension of a finite group by a nilpotent group. In particular, all nilpotent groups are finite-by-nilpotent (taking $K= \{e\}$). Moreover, the class of finite-by-nilpotent groups is strictly larger than the class of nilpotent groups, as every finite group belongs to the former class but not to the latter.

 Finite-by-nilpotent groups also admit a nice description in terms of their lower central series: it is known that finite-by-nilpotent groups are precisely those whose lower central series stabilizes at a finite group.

We are now going to show that for any Hecke pair $(G,\gm)$ where $G$ is finite-by-nilpotent, every double coset $\gm s \gm$ generates a finite co-hereditary set, implying that $C^*(G, \gm)$ exists and coincides with $C^*(L^1(G, \gm))$.

 Let $s \in G$ and $\{\gamma_k\}_{k \in \mathbb{N}} \subseteq \gm$. It is clear that $[s, \gamma_1, \dots, \gamma_k] \in G_k $. Since the series $\{G_k\}$ eventually stabilizes at a finite subgroup, it follows directly from Corollary \ref{D_k finito} $a)$ that $\gm s \gm$ generates a finite co-hereditary set. This concludes the proof.\\

\subsection{$G$ is Hypercentral}

Recall that a group $G$ is said to be a \emph{hypercentral group} (also called a \emph{$ZA$-group}) if its upper central series, possibly continued transfinitely, stabilizes at the whole group $G$. For a rigorous definition of this concept, we refer the reader to \cite[section 12.2]{rob2} for example. Another characterization of hypercentral groups, which is the one we will use, is given by the following result:\\

\begin{thm}[Lemma, page 219, \S 63, \cite{kurosh}]
\label{hypercentral groups kurosh}
 A group $G$ is hypercentral if and only if it satisfies the following property: for any $s \in G$ and any sequence $\{x_n\}_{n \in \mathbb{N}} \subset G$ there is a $k \in \mathbb{N}$ such that
\begin{align*}
 [s,x_1, \dots, x_k]=e\,.\\
\end{align*}
\end{thm}

We will now prove that if $(G,\gm)$ is a Hecke pair with $G$ a hypercentral group, then every double coset $\gm s \gm$ generates a finite co-hereditary set, so that $C^*(G, \gm)$ exists and $C^*(G, \gm) \cong L^1(G, \gm)$. This is a direct application of Corollary \ref{D_k finito} $a)$, taking $F = \{e\}$, given the characterization of hypercentral groups of Theorem \ref{hypercentral groups kurosh}.\\

\begin{rem}
\label{remark hypercentral}
 The class of hypercentral groups also strictly contains the class of nilpotent groups (see Example \ref{2 -quasicyclic group}), and moreover it is known that every hypercentral group is locally nilpotent (but not vice-versa). Thus, we have found another class of groups $G$, satisfying a nilpotent-type property, for which the Hecke algebra $\h(G, \gm)$ of any Hecke pair $(G, \gm)$ has an enveloping $C^*$-algebra (which coincides with $C^*(L^1(G, \gm))$).\\
\end{rem}

\begin{ex}
\label{2 -quasicyclic group}
 Let $\mathbb{Z}_{2^{\infty}}$ be the $2$-quasicyclic group, i.e. the group of all the $2^n$-th roots of unity for all $n \in \mathbb{N}$. This group is the Pontryagin dual of group of $2$-adic integers. The group $\mathbb{Z} / 2 \mathbb{Z}$ acts on $\mathbb{Z}_{2^{\infty}}$ by mapping an element to its inverse. The generalized dihedral group
\begin{align*}
G := \mathbb{Z}_{2^{\infty}}\rtimes \big(\mathbb{Z} / 2 \mathbb{Z} \big)
\end{align*}
is a group which is hypercentral, but not nilpotent.\\
\end{ex}

\subsection{$G$ is an $FC$-group and $\gm$ is Finite}

Recall that a group $G$ is said to be $FC$ if every element $s$ has finitely many conjugates, i.e. the set $\mathcal{C}_s :=\{t^{-1}s t: t \in G \}$ is finite. It can be seen that every subgroup $\gm \subseteq G$ of an $FC$-group is a Hecke subgroup, because
\begin{align*}
 \gm s \gm = \bigcup_{\gamma \in \gm} \gamma s \gm = \bigcup_{\gamma \in \gm} \gamma s \gamma^{-1} \gm \subseteq \bigcup_{x \in \mathcal{C}_s} x \gm\,,
\end{align*}
and the last union is finite.

$FC$ groups are a generalization of both finite and abelian groups, and share many common properties with these classes. They were extensively studied by B. H. Neumann and others, starting with the article \cite{neumann}. The analogous class of groups in the locally compact setting (groups in which the conjugacy class of any element has compact closure) is usually denoted by $FC^{-}$ and has also been widely studied, since it is a direct generalization of both compact and abelian locally compact groups (see \cite[Chapter 12]{palmer} for an account).

When $G$ is a $FC$-group and $\gm \subseteq G$ is a finite subgroup, we can prove that every double coset $\gm s \gm$ generates a finite co-hereditary set, so that $C^*(G, \gm)$ exists and $C^*(G, \gm) \cong C^*(L^1(G, \gm))$. To see this, let $s \in G$ and $\{\gamma_k \}_{k \in \mathbb{N}} \subseteq \gm$. Also, let $\gm = \{\theta_1, \dots, \theta_n\}$ and for each $1 \leq i \leq n$ let us denote by $S_i \subseteq \mathbb{N}$ the set
\begin{align*}
 S_i:=\{j \in \mathbb{N} : \gamma_j = \theta_i\}\,.
\end{align*}
Of course, the sets $S_i$ are mutually disjoint and their union is $\mathbb{N}$. We have that
\begin{eqnarray*}
 \{\gm [s, \gamma_1, \dots, \gamma_k] \gm: k \in \mathbb{N}\} & = & \bigcup_{i =1}^n \{\gm [s, \gamma_1, \dots, \gamma_k] \gm: k \in S_i\}\\
& = & \bigcup_{i =1}^n \{\gm [s, \gamma_1, \dots, \gamma_{k-1}, \theta_i] \gm: k \in S_i\}\,.
\end{eqnarray*}
Now we notice that
\begin{align*}
 \gm [s, \gamma_1, \dots, \gamma_{k-1}, \theta_i]\gm = \gm [s, \gamma_1, \dots, \gamma_{k-1}]^{-1} \theta_i^{-1} [s, \gamma_1, \dots, \gamma_{k-1}]\gm\,.
\end{align*}
 Since there are only finitely many conjugates of $\theta_i^{-1}$, it follows that the set $\{\gm [s, \gamma_1, \dots, \gamma_{k-1}, \theta_i] \gm: k \in S_i\}$ is finite, and therefore $\{\gm [s, \gamma_1, \dots, \gamma_k] \gm: k \in \mathbb{N}\}$ is finite. Thus, by Theorem \ref{grafo e serie}, the co-hereditary set generated by $\gm s \gm$ is finite.\\

\subsection{$G$ is Locally-Nilpotent and $\gm$ is Finite}

Recall that a group $G$ is said to be \emph{locally-nilpotent} if every finitely generated subgroup of $G$ is nilpotent.

Let $G$ be a locally-nilpotent group and $\gm$ a finite subgroup. The pair $(G,\gm)$ is automatically a Hecke pair since $\gm$ is finite. We are now going to prove that each double coset $\gm s \gm$ generates a finite co-hereditary set, implying that $C^*(G, \gm)$ exists and coincides with $C^*(L^1(G, \gm))$. To see this, let $\langle s, \gm \rangle \subseteq G$ be the subgroup generated by $s$ and $\gm$. This subgroup is finitely generated, hence nilpotent. Thus, as we have proven above, $\gm s \gm \in \h(\langle s , \gm \rangle, \gm) \subseteq \h(G, \gm)$ generates a finite co-hereditary set.\\

\subsection{$G$ is Locally-Finite and $\gm$ is Finite}
\label{loc finite classes}

Recall that a group $G$ is said to be \emph{locally-finite} if every finitely generated subgroup of $G$ is finite. 

Let $G$ be a locally-finite group and $\gm$ a finite subgroup. The pair $(G,\gm)$ is automatically a Hecke pair since $\gm$ is finite. We are now going to prove that each double coset $\gm s \gm$ generates a finite co-hereditary set, implying that $C^*(G, \gm)$ exists and coincides with $C^*(L^1(G, \gm))$. To see this, let $\langle s, \gm \rangle \subseteq G$ be the subgroup generated by $s$ and $\gm$. This subgroup is finitely generated, hence finite. Thus, as we have proven above, $\gm s \gm \in \h(\langle s , \gm \rangle, \gm) \subseteq \h(G, \gm)$ generates a finite co-hereditary set.\\

An interesting feature of Hecke pairs arising from locally finite groups is that they give rise to AF Hecke algebras. In that regard we have the following result:\\

\begin{prop}
Let $(G, \gm)$ be a Hecke pair where $G$ is countable and $\gm$ is a finite subgroup. Then $\h(G, \gm)$ is an AF $^*$-algebra if and only if $G$ is locally finite.\\
\end{prop}

{\bf \emph{Proof:}} $(\Longleftarrow)$ Assume $G$ is locally finite. Since $G$ is assumed countable, let us fix an enumeration of its elements $G =\{g_1, g_2, \dots\}$ and for each $n \in \mathbb{N}$ let us define $H_n$ as the subgroup $H_n := \langle \gm, g_1, \dots, g_n \rangle$. It is clear that $\{H_n\}_{n \in \mathbb{N}}$ forms an increasing sequence of finitely generated subgroups, such that $\bigcup H_n = G$. Moreover, since $G$ is locally finite, each $H_n$ is a finite group which contains $\gm$. Hence, we have a sequence of finite dimensional Hecke algebras $\{\h(H_n, \gm)\}_{n \in \mathbb{N}} \subseteq \h(G, \gm)$ satisfying $\bigcup \h(H_n, \gm) = \h(G, \gm)$. Thus, $\h(G, \gm)$ is an $AF$ $^*$-algebra.

$(\Longrightarrow)$ Assume that $\h(G, \gm)$ is an AF $^*$-algebra. Then any element $f \in \h(G, \gm)$ lies in a finite dimensional $^*$-subalgebra, and is therefore algebraic over $\mathbb{C}$. It then follows from \cite[Proposition 2.6]{krieg} that $G$ is locally finite. \qed\\

\begin{ex}
 Similarly to Example \ref{2 -quasicyclic group}, let $p$ be a prime number and $\mathbb{Z}_{p^{\infty}}$ be the $p$-quasicyclic group (which is the Pontryagin dual of the group of $p$-adic integers). The generalized dihedral group
\begin{align*}
G := \mathbb{Z}_{p^{\infty}}\rtimes \big(\mathbb{Z} / 2 \mathbb{Z} \big)\,,
\end{align*}
is locally finite (but not locally nilpotent unless $p =2$).\\
\end{ex}

\section{Reduction Techniques}
\label{reduction section}

Suppose we have $\gm \subseteq N$ for some normal subgroup $N$ of $G$, i.e.
\begin{align*}
 \gm \subseteq N \unlhd G\,.
\end{align*}
We will see that we can ``reduce'' the problem of knowing whether $\mathcal{H}(G,\gm)$ has an enveloping $C^*$-algebra to the smaller Hecke $^*$-algebra $\mathcal{H}(N,\gm)$. This kind of problem has been considered in the literature (see \cite[Proposition 2.6]{compac} and also \cite[Theorem 1.11]{laca larsen} with a more general version given in \cite[Theorem 5.5]{grp ext}). Sometimes, as presented in the last two references, the structure of $\mathcal{H}(G,\gm)$ can be described in terms of $\mathcal{H}(N,\gm)$. At this point, however, we will address solely the question of existence of enveloping $C^*$-algebras, for which our approach is more general. We have the following result:\\

\begin{prop} Let $(G,\gm)$ be a Hecke pair and $N$ a normal subgroup of $G$ containing $\gm$. If $\mathcal{H}(N,\gm)$ has an enveloping $C^*$-algebra (resp. is a $BG^*$-algebra), then the same is true for $\mathcal{H}(G,\gm)$. Moreover, all double cosets of $\mathcal{H}(N,\gm)$ generate finite co-hereditary sets if and only if all double cosets of $\mathcal{H}(G,\gm)$ do so.\\
\label{(N,H) (G,H)}
\end{prop}

{\bf \emph{Proof:}} Since $\mathcal{H}(N,\gm)$ is a $^*$-subalgebra of $\mathcal{H}(G,\gm)$, all $^*$-representations of $\mathcal{H}(G,\gm)$ restrict to $^*$-representations of $\mathcal{H}(N,\gm)$ (and similarly for pre-$^*$-representations).

Suppose that $\mathcal{H}(N,\gm)$ has an enveloping $C^*$-algebra. Let $\gm s\gm \in \mathcal{H}(G,\gm)$. By Remark \ref{g gamma h in prod of double cosets} we know that
\begin{align*}
 (\gm s \gm)^**\gm s\gm = \sum_{i=1}^n \lambda_i\, \gm s^{-1} \gamma_i s\gm
\end{align*}
for some $n \in \mathbb{N}$, $\lambda_1, \dots, \lambda_n \in \mathbb{N}$ and $\gamma_1, \dots, \gamma_n \in \gm$. If $\pi$ is a $^*$-representation of $\mathcal{H}(G,\gm)$, we have that
\begin{eqnarray*}
 \|\pi(\gm s\gm)\| & = & \sqrt{\|\pi((\gm s\gm)^**\gm s\gm)\|} \;\;\leq\;\; \sqrt{\sum_{i=1}^n \lambda_i \|\pi(\gm s^{-1} \gamma_i^{-1} s\gm)\|}\,.
\end{eqnarray*}
 We notice that, since $N$ is normal in $G$ and $\gm \subseteq N$, we have $s^{-1} \gamma_i^{-1} s \in N$, i.e. all elements $\gm s^{-1} \gamma_i^{-1}s\gm \in \mathcal{H}(N,\gm)$. Also, by restriction, $\pi$ gives rise to a $^*$-representation of $\mathcal{H}(N,\gm)$. Hence, we must have
\begin{align}
 \|\pi(\gm s\gm)\| \leq \sqrt{\sum_{i=1}^n \lambda_i \|\gm s^{-1}\gamma_i^{-1}s\gm\|_{u, N}}\,,
\label{desigualdade}
\end{align}
where $\| \cdot \|_{u, N}$ is the universal norm in $\mathcal{H}(N,\gm)$. Since the inequality in (\ref{desigualdade}) holds for every $^*$-representation $\pi$, we see that $\gm s\gm$ has a bounded universal norm. Thus, $\mathcal{H}(G,\gm)$ has an enveloping $C^*$-algebra.

A similar computation as the above, but with pre-$^*$-representations instead, shows that if $\h(N, \gm)$ is a $BG^*$-algebra, then so is $\h(G, \gm)$.

Let us now prove the second statement. First we notice that if all double cosets of $\mathcal{H}(G,\gm)$ generate finite co-hereditary sets, then the same is true for $\mathcal{H}(N,\gm)$, since $\mathcal{H}(N,\gm) \subseteq \mathcal{H}(G,\gm)$. Now suppose that all double cosets of $\mathcal{H}(N,\gm)$ generate finite co-hereditary sets. Then, given any $\gm s \gm \in \mathcal{H}(G,\gm)$, we have that all of its successors are of the form $\gm s^{-1} \gamma s\gm$ for some $\gamma \in \gm$, and therefore belong to $\mathcal{H}(N,\gm)$, as we saw above. Notice that the co-hereditary set generated by $\gm s\gm$ is the union of the co-hereditary sets generated by its successors (along with the element $\gm s \gm$ itself). Thus, since its successors belong to $\h(N, \gm)$, they generate finite co-hereditary sets, and since $\gm s \gm$ has finitely many successors, $\gm s \gm$ also generates a finite co-hereditary set. \qed\\

An easy induction argument allows us to extend the previous result to ascendant subgroups, and thus, in particular, to subnormal subgroups (recall the definitions given in the sections \ref{ascendant section} and \ref{subnormal section}): \\

\begin{cor}
 Let $(G,\gm)$ be a Hecke pair and $K$ an ascendant subgroup of $G$ containing $\gm$.  If $\mathcal{H}(K,\gm)$ has an enveloping $C^*$-algebra (resp. is a $BG^*$-algebra), then the same is true for $\mathcal{H}(G,\gm)$. Moreover, all double cosets of $\mathcal{H}(K,\gm)$ generate finite co-hereditary sets if and only if all double cosets of $\mathcal{H}(G,\gm)$ do so.\\
\label{(K,H) (G,H) subnormal}
\end{cor}

By combining Proposition \ref{(N,H) (G,H)}, and the more general Corollary \ref{(K,H) (G,H) subnormal}, together with the results of the previous section, we obtain many more classes of Hecke pairs for which the existence of an enveloping $C^*$-algebra can be assured. For example, if $(G, \gm)$ is a Hecke pair for which there exists an ascendant nilpotent subgroup $K$ with $\gm \subseteq K \subseteq G$, then $\h(G, \gm)$ has an enveloping $C^*$-algebra. Of course, one can replace ``nilpotent'' with any other property discussed in the previous section.

As yet another example, we can recover (and slightly improve) a result by Baumgartner, Ramagge and Willis (\cite[Proposition 2.6]{compac}), concerning Hecke pairs $(G, \gm)$ where $\gm$ has finite index in a normal subgroup $N\unlhd\, G$. They proved that in this case $\h(G, \gm)$ has an enveloping $C^*$-algebra. Through our Corollary \ref{(K,H) (G,H) subnormal} we can improve this result in the following way:\\

\begin{cor}
Let $(G, \gm)$ be a Hecke pair such that $\gm$ has finite index in an ascendant subgroup $K \subseteq G$. Then $\h(G, \gm)$ has an enveloping $C^*$-algebra which coincides with $C^*(L^1(G, \gm))$.\\
\end{cor}

\section{Final Remarks and Questions}
\label{questions section}

\begin{enumerate}
 \item In this article we studied conditions that imply the existence of a full Hecke $C^*$-algebra. What about the opposite question: when can we ensure that $C^*(G, \gm)$ does not exist for a given Hecke pair? Even though there are examples of Hecke pairs $(G, \gm)$ for which it is known that $C^*(G, \gm)$ does not exist (see \cite[Proposition 2.21]{hall} or \cite[Example 3.4]{tzanev}), this question seems to be quite delicate in full generality. It would be interesting to have sufficient (group-theoretic) conditions for the non-existence of $C^*(G, \gm)$.

\item If $C^*(G, \gm)$ exists, does it necessarily follow that $C^*(G, \gm) \cong C^*(L^1(G, \gm))$? With the notable exception of the Iwahori Hecke algebras (section \ref{iwahori hecke class}), for which we do not know the answer, this is true for all the other examples we know of. 

\item An even stronger question, posed in \cite{schl}, is still open: if $C^*(G, \gm)$ exists, does it necessarily follow that $C^*(G, \gm) \cong C^*(L^1(G, \gm)) \cong pC^*(\overline {G})p$ ? Many of the new classes studied in this article will be shown in a following article (see \cite{palma2}) to have this property.

\item Is it true that if $C^*(G, \gm)$ exists, then $\h(G, \gm)$ is a $BG^*$-algebra? All examples of Hecke algebras for which we know that an enveloping $C^*$-algebra exists are actually $BG^*$-algebras. It would be interesting to have a counter-example or a proof of this statement.
\end{enumerate}

\end{document}